\documentclass[a4paper,12pt]{amsart}

\usepackage[all]{xy}
\usepackage{amscd}
\usepackage{amsfonts}
\usepackage{a4wide}
\usepackage{amsmath}
\usepackage{amssymb}
\usepackage{natbib}
\usepackage{fullpage}
\usepackage[applemac]{inputenc}

\usepackage[pdftex,backref,hyperindex,breaklinks,colorlinks,citecolor=blue]{hyperref}

\numberwithin{equation}{section}

\newcommand{\sg}{\sigma}
\newcommand{\eps}{\varepsilon}
\newcommand{\ph}{\varphi}

\newcommand{\Fp}{\mathbb{F}_p}

\newcommand{\ce}{\mathbb{C}}
\newcommand{\erre}{\mathbb{R}}
\newcommand{\z}{\mathbb{Z}}
\newcommand{\ene}{\mathbb{N}}
\newcommand{\q}{\mathbb{Q}}

\newcommand{\h}{\mathbb{H}}
\newcommand{\p}{\mathbb{P}^1}
\newcommand{\disco}{\mathbb{D}}

\newcommand{\XC}[1]{X_0 (#1)}
\newcommand{\XQ}[1]{X_0 (#1)_\q}

\newcommand{\XZ}[1]{\mathcal{X}_0 (#1)}
\newcommand{\YC}[1]{Y_0 (#1)}

\newcommand{\G}{\Gamma_0}

\newcommand{\Ch}{\widehat{CH} (N)}
\newcommand{\ChR}{\widehat{CH}_\erre (N)}
\newcommand{\Chnum}{\widehat{CH}_\erre^{num} (N)}

\newcommand{\infinito}{\hat{D}_\infty}
\newcommand{\JNQ}{J_0 (N)}

\newcommand{\JN}{J_0 (N)}

\newcommand{\Eis}{\textrm{Eis}}

\newcommand{\X}{\mathcal{X}}
\newcommand{\Y}{\mathcal{Y}}
\newcommand{\CH}{\widehat{CH^1}}
\newcommand{\dvemph}{\textrm{\emph{div} }}
\newcommand{\dv}{\textrm{div }}

\newcommand{\der}{\mathclose>}
\newcommand{\izq}{\mathopen<}

\newcommand{\Sch}{\mathcal{D}}
\newcommand{\wbar}{\hat{\omega}}

\newcommand{\dz}{dz \wedge d \bar{z}}

\newcommand{\Luno}{L^2_1}

\newcommand{\pruf}{\noindent \textbf{Proof}: }
\newcommand{\esp}{L^2_{-1}}
\newcommand{\espa}{\mathcal{L}}

\newtheorem{thm}{Theorem}[section]
\newtheorem{lema}[thm]{Lemma}
\newtheorem{prop}[thm]{Proposition}
\newtheorem{prop-def}[thm]{Proposition-Definition}
\newtheorem{cor}[thm]{Corollary}
\theoremstyle{definition}
\newtheorem{rem}[thm]{Remark}
\newtheorem{rems}[thm]{Remarks}
\newtheorem{ex}[thm]{Example}
\newtheorem{definicion}[thm]{Definition}
\newtheorem{definiciones}[thm]{Definitions}

\title{Correspondences in  Arakelov geometry and applications to the case of Hecke operators on modular curves}
\author{Ricardo Menares}
\address{Facultad de Matem\'aticas \\
Pontificia Universidad Cat\'olica de Chile \\
Vicuna Mackenna 4860 \\
Santiago, Chile}
\email{rmenares@mat.puc.cl}

\begin{document}

\begin{abstract}
In the context of arithmetic surfaces,  Bost defined a generalized Arithmetic Chow Group
 (ACG) using the Sobolev space $\Luno$. We study the behavior of these groups under pull-back and push-forward and we prove a projection formula.
 We use these results to define an action of the Hecke operators on the ACG of modular curves  and  to show that they are self-adjoint with respect to the arithmetic 
intersection product.    The decomposition of the ACG in eigencomponents which follows  allows us to define new numerical invariants, which are refined versions of 
 the self-intersection of the dualizing sheaf. Using  the Gross-Zagier formula and a calculation due independently to Bost and K\"uhn we compute these invariants in terms 
of special values of L series. On the other hand, we obtain a proof of the fact that Hecke correspondences acting on the Jacobian of the modular curves are self-adjoint with
respect to the N\'eron-Tate height pairing.
\end{abstract}

\maketitle
\tableofcontents

\section{Introduction}

Let $\X$ be an arithmetic surface, that is a  projective, regular, 2-Krull dimensional scheme over Spec$(\z)$. Such an object arises naturally  as a model over the integers of a smooth algebraic curve defined over a number field.  Arakelov devised  a formalism on $\X$ analogous to the classical intersection theory on an algebraic surface defined over a field (\cite{Arakelov}). The basic object of this formalism is a compactified divisor, that is  a pair $(D,g)$, where $D$ is a Weil divisor on $\X$ and $g$ is a real-valued distribution on the compact Riemann surface $X=\X(\ce)$. Here, $g$ is a Green function for $D$, subject to the following condition. Suppose for simplicity that the genus of  $X$ is positive, and take an orthonormal basis  $\omega_1, \ldots, \omega_n$ of the space of holomorphic 1-forms  (where we take this space to be endowed with the inner product $(\alpha, \beta) \mapsto \frac{i}{2}\int_{X}^{}\alpha \wedge \overline{\beta}$). The  canonical Arakelov 2-form is defined by $$\mu_{Ar}:=\frac{i}{n}\sum_{j=1}^{n}\omega_j \wedge \overline{\omega_j}.$$ This is a $C^\infty$ volume form on $X$. It is independent of the choice of the orthonormal basis $\{\omega_1, \ldots, \omega_n\}$.  We impose that  $g$ satisfies the equation of currents 

\begin{equation}\label{defArakelov}
dd^c g + \delta_{D(\ce)} = \big(\deg D(\ce)\big) \mu_{Ar}.
\end{equation}

\noindent If this condition is satisfied, then $g$ is said to be admissible. We also impose the normalization $\int_{X}^{} g \mu_{Ar}=0$.  

The key analytic ingredient needed to define an intersection number from two compactified divisors $(D_1,g_1)$ and $(D_2,g_2)$ is the star product 

\begin{equation} \label{star}
g_1 * g_2 =  g_1(dd^c g_2 + \delta_{D_2(\ce)})  +g_2 \delta_{D_1(\ce)}.
\end{equation}

\noindent This is a well defined current if $D_1(\ce)$ and $D_2(\ce)$ do not have points in common.
The integral $\int_{\X(\ce)}^{}g_1 * g_2$ is bilinear and symmetric in $(g_1,g_2)$.

In  Gillet and Soulé's subsequent treatment of Arakelov's theory \cite{GilletSoule}, the admissibility condition on $g$ is dropped. More precisely, $g$ is a real-valued $C^\infty$ function on the complement of the support of $D(\ce)= \sum_{P}^{}a_P [P]$. It is required that for every point $P$ in this support,  there is a local expansion of the form 

\begin{equation}\label{defGilletSoule}
g = -a_P\log |\theta_P(\cdot)|^2 + b,
\end{equation}
\noindent  with $\theta_P$  a local chart vanishing at $P$ and $b$  a $C^\infty$ function. Such a distribution is called a \emph{Green current} for $D$. Equation \eqref{defArakelov} implies the local expansion \eqref{defGilletSoule}, making Gillet and Soulé's theory a generalization of Arakelov's theory. The star product \eqref{star} between normalized Green currents  retains the symmetry and bilinearity properties.

In this paper we study the  \emph{lack of functoriality} problem in Arakelov theory. That is, let $\ph : \X \rightarrow \Y$ be a finite morphism between arithmetic surfaces, and $(D,g)$ a compactified divisor on $\X$. If the induced  holomorphic covering $\ph_\ce : \X(\ce) \rightarrow \Y(\ce)$ is ramified, then the pair  $(\ph_* D,  \ph_* g)$ induced  from the the push-forward is not a compactified divisor on $\Y$, in either the Arakelov or the Gillet-Soulé sense. The reason is  purely analytic in nature; $\ph_* \mu_{Ar}$ has singularities at the branching points, hence it is not the Arakelov canonical form on $\Y(\ce)$. Furthermore,  $\ph_* g$ is not a Green current for $\ph_*D$, because the push-forward operation does not preserve $C^\infty$ functions. In this sense, the theories of Arakelov and Gillet-Soulé are not functorial. In many situations, it is natural to consider correspondences on arithmetic surfaces. Unfortunately, the interaction with the intersection formalism is excluded by this phenomenon.

Bost brought to this circle of ideas   the Sobolev space $\Luno$, which consists of distributions in $L^2$ whose derivatives  also belong to $L^2$. He proposed a generalization of Gillet-Soulé's theory for arithmetic surfaces  (\cite{Bost}). In this framework, $g$ is a real-valued  element of $\Luno$ away from the support of $D(\ce)$ and for every point $P$ in this support, there is a local expansion similar to \eqref{defGilletSoule} but where $b$ is in $\Luno$. Such an objet is called a \emph{$\Luno$-Green function}.  The integral of the star product \eqref{star} between $\Luno$-Green functions is well defined. This theory is optimal in the sense that $\Luno$ is the biggest space such that a useful (integral of a) star product can be defined. Bost's theory also improves on the geometrical part, because it doest not require the scheme $\X$ to be regular but only integral and normal. We drop the regularity in the definition of arithmetic surface in what follows.

Let us write  $\big(\CH(\X), \izq, \der \big)$ to denote the $\Luno$ arithmetic Chow group defined in \cite{Bost}. In this paper, we prove that these groups have a good functoriality. More precisely, we provide a proof of the following

\begin{thm}\label{main intro}
Let $\ph: \X \rightarrow \Y$ be a generically finite morphism between integral, normal, arithmetic surfaces.
\begin{enumerate}

\item \label{bien definido intro} The natural definitions of pull-back and push-forward induce homomorphisms 
\begin{eqnarray}
\varphi^*  :   \CH(\Y) \longrightarrow \CH(\X)  \nonumber \\
\varphi_*  :   \CH(\X) \longrightarrow \CH(\Y).  \nonumber 
\end{eqnarray}  

\item \label{formula de proyeccion intro} Projection formula: for all $x \in \CH (\X)$ and $y \in \CH(\Y)$, we have that
\[
\izq\varphi^* y, x\der = \izq  y, \varphi_*x\der.
\]
\end{enumerate}

\end{thm}

The main difficulty in proving this theorem is to control the $L^2$ norm of $\ph^* f$ around the ramification points, where $f$ is $\Luno$. This is achieved by establishing a Poincaré-like inequality (Lemma \ref{ramificacion} in the main text), which we present as a variant of the Hardy inequality. On the other hand, the push-forward operation turns out to be easier to handle. 

We point out that Theorem \ref{main intro} was known to Bost before our investigation started (\cite{Bost}, p. 245, \cite{carta}), though  he did not publish any details. The method described in the paragraph above, different from Bost's approach, was found independently (cf. Remark \ref{usando DL}).

We also consider a generalization of Arakelov's condition \eqref{defArakelov}. Let $\esp \big(\X(\ce)\big)$ be the space of 2-currents which are locally of the form  $(f + \frac{\partial}{\partial z} h_1 + \frac{\partial}{\partial \bar{z}} h_2)\dz$,  where $f, h_1, h_2$ are $L^2$. Let $\mu$ be a real-valued element of  $\esp \big(\X(\ce)\big) $ such that $\int_{\X(\ce)}^{}\mu=1$. Given a Weil divisor $D$ on $\X$, we call a $\mu$-admissible Green function for $D$ an 
$\Luno$-Green function $g$ for $D$ such that  equation $$dd^c g + \delta_{D(\ce)} = \big( \deg D(\ce)\big) \mu $$ holds in the sense of currents.  We denote by $\CH(\X)_\mu$ the corresponding $\mu$-admissible arithmetic Chow group.
We show that pull-back and push-forward are well-behaved with respect to this notion, as well as study how the reference 2-current changes. 

\begin{thm} \label{funtorial intro} Let $\mu \in \esp \big(\X(\ce)\big)$ and  $\nu \in \esp \big(\Y(\ce)\big)$ be real-valued 2-currents such that $\int_{\X(\ce)}^{}\mu = \int_{\Y(\ce)}^{}\nu=1$.  The natural definitions of pull-back and push-forward induce homomorphisms

\begin{eqnarray}
\varphi^*  :   \CH(\Y)_\nu \longrightarrow \CH(\X)_{\frac{\varphi^*\nu}{\deg \ph}}  \nonumber \\
\varphi_*  :   \CH(\X)_\mu \longrightarrow \CH(\Y)_{\varphi_*\mu}.  \nonumber 
\end{eqnarray}  
\end{thm}

Here, the main issue is the question of the existence of such $\mu$-admissible Green functions (cf. Proposition \ref{existencia}). Such a distibution always exists locally, a result coming  from the theory of elliptic PDEs.  This reduces our problem to patching together local solutions to obtain a global one, which in turn is an exercise in $\Luno$-cohomology of Riemann surfaces. This exercise is probably well known to the experts, though for lack of suitable reference in the literature, we collect in the appendix the necessary facts and deductions to complete it.

We mention that Burgos, Kramer and Kühn \cite{BKK} have developed an arithmetic intersection theory which, in the case of arithmetic surfaces considered here, lies strictly in between Gillet-Soulé's and Bost's. They show that their theory is functorial. This theory has the advantage of being able to handle higher dimensional arithmetic varieties, where no optimal theory analogous to Bost's is known (though there are partial results in this direction due to A. Moriwaki \cite{Moriwaki}).

As an application, which is also our main motivation to study functoriality in this context, we consider the case of the modular curves $\XC{N} = \G(N) \backslash \h^*$ and their Hecke correspondences.  We use the  integral model $\XZ{N}$ obtained by the modular interpretation in terms of $\G(N)$-structures  as in \cite{KatzMazur}. On the other hand, the hyperbolic measure on $\h$ descends on $\XC{N}$ to a finite measure. We denote by $\mu_N$ the normalized 2-current on $\XC{N}$. This is an element of $\esp \big(\XC{N}\big)$. We define $\Ch := \CH (\XZ{N})_{\mu_N}$.

\begin{thm}
The Hecke correspondences $T_l$ with $l \nmid N$ and the involutions $w_d$ with $d | N$ induce homomorphisms $\hat{T}_l$ and $\hat{w}_d$ on $\CH (\XZ{N})$. These homomorphism form a commutative algebra, and the induced operator $\hat{T}_l$ is self-adjoint. Moreover, this  algebra preserves $\Ch$.
\end{thm}

We remark that, as the singularities of an admissible Green function with respect to $\mu_N$ are of $\log\log$ type, the part of this  statement concerning functoriality can also be obtained via the Burgos-Kramer-Kühn theory.

From now on, we suppose that $N$ is squarefree.  We denote by $\Chnum$ the variant of the arithmetic Chow group $\Ch$ obtained by taking divisors (compactified, $\mu_N$-admissible) with real coefficients, and then the quotient by the numerical equivalence relation (cf. Section \ref{num} for a precise definition). This is a real vector space of finite dimension, and the Hecke operators $\hat{T}_l$ and $\hat{w}_d$ act on it. The Hodge index theorem in this context allows us to decompose $\Chnum$ as a sum of two Hecke-invariant subspaces:

\begin{equation} \label{primeradescomposicion}
\Chnum = \Eis \overset{\perp}{\oplus} J,
\end{equation} 
\noindent where $\Eis$ is the space spanned by compactified irreducible components of fibers, and 

\begin{equation} \label{jotas}
J \cong \JN(\q)\otimes \erre
\end{equation} 

\noindent as real vector spaces. It is possible to choose the isomorphism in \eqref{jotas} such that the Hecke actions on both sides are compatible (Corollary \ref{compatibilidad}). Hence we may further decompose

\begin{equation} \label{descomposicionintro}
J = \bigoplus_f J_f, 
\end{equation}

\noindent where $f$ runs through a basis of $S_2\big(\G (N)\big)$ consisting of eigenforms.

Let $\omega$ be the dualising sheaf on $\XZ{N}$. By taking the divisor induced by a section of this sheaf with an appropriate choice of $\Luno$-Green function (cf. Section \ref{explicacion de la normalizacion}), we obtain a class  $\wbar \in \Ch$. Following the decompositions \eqref{primeradescomposicion} and \eqref{descomposicionintro} we can then write

\[
\hat{\omega} =\hat{\omega}_{\Eis} +  \hat{\omega}_{J}, \quad \hat{\omega}_{J} = \sum_{f}^{} \hat{\omega}_{f}.
 \]

We consider the self-intersection of each component as a refined invariant: 

\[\hat{\omega}^2 = \hat{\omega}_{\Eis}^2 + \hat{\omega}_{J}^2, \quad \hat{\omega}_{J}^2 = \sum_{f}^{} \hat{\omega}_{f}^2.\]

Using a calculation due to U. Kühn (\cite{Kuhn}), we then obtain

\begin{thm} \label{calculo Eis} Suppose that $N$ is squarefree. Let $g_M := \textrm{genus of  } \XC{M}$. We have that

\begin{eqnarray}
\frac{\wbar_{\emph{\Eis}}^2}{(g_N-1)^2} = \frac{576}{[\G(1) : \G(N)]} \Big( \frac{1}{2} \zeta(-1)
 + \zeta'(-1)\Big) - \sum_{p|N} \frac{\log p}{g_N-2 g_{\frac{N}{p}} +1} ,
\end{eqnarray}

\noindent where $\zeta(\cdot)$ is the Riemann zeta function.

\end{thm}

One of the  divisors in the class $\omega_J$ is an explicit  linear combination of Heegner points of discriminants -3 and -4. Let $H_i$ (resp. $H_j$) be the sum of all divisors on $\XC{N}$ of the form $\frac{1}{2}\big([P]-[\infty]\big)$ (resp. $\frac{1}{3} \big( [P]-[\infty]\big)$),  with $P$ a Heegener point of discriminant -4 (resp. -3). Using the Gross-Zagier formula, we obtain that

\begin{thm} \label{calculo isotipico} Suppose that $f$ is a normalized newform, and that $N$ is a squarefree integer such that $(6,N)=1$. Then,

\begin{equation} \label{alturas}
 \wbar_f^2 = - \Big( h_{NT}\big((H_i)_f\big)^{1/2} +  2 h_{NT}\big((H_j)_f\big)^{1/2} \Big)^2.
\end{equation}

\noindent where   $h_{NT}$ is the Néron-Tate height on $\JN(\q)$. Furthermore, if $N$ is prime, then

$$h_{NT} (H_i)_f = \frac{1}{2\pi^2(f,f)} L(f, \chi_{-4},1) L'(f,1), \quad h_{NT} (H_j)_f = \frac{\sqrt{3}}{4\pi^2(f,f)}  L(f, \chi_{-3},1) L'(f,1),$$ where $\chi_{-4}$ (resp. $\chi_{-3}$) is the quadratic Dirichlet character associated to $\q(i)/\q$ (resp. $\q(\sqrt{-3})/\q$) and $(f,f)$ is the Petersson norm of $f$.

\end{thm}

Though we have treated in detail only the case of $\G(N)$,  an analogous analysis can be made for Hecke operators on other  congruence groups, as well as on Shimura curves.

\noindent \textbf{Acknowledgements}: This work is  based in part on the author's Ph.D. thesis at Université d'Orsay, written  under the supervision of E.  Ullmo. It was A. Chambert-Loir who suggested  viewing the results of the thesis as a functoriality problem on an appropiate function space. We have benefited from  useful concrete suggestions by J.-B. Bost, U. Kühn and A. Thuillier.

\section{Functoriality in $\Luno$}

\subsection{The space $\Luno$ and related spaces}
We recall the basic definitions from \cite{Bost}.
Let $X$ be a compact Riemann surface. The space $\Luno (X)$ is the space of distributions  $f \in L^2 (X)$ such that $\frac{\partial}{\partial z} f$ and $\frac{\partial}{\partial \bar{z}} f$ are also in $L^2(X)$.   We have a hermitian inner product $(\cdot, \cdot)_{1,X}$ on $\Luno (X)$ given by 

\[
(f, g)_{1,X} = i  \int_{X}^{} \partial f \wedge \overline{\partial g}.
\]

\noindent Its associated seminorm is given by

\[
\|f\|_{1,X}^2 = i  \int_{X}^{} \partial f \wedge \overline{\partial f} \in \erre_{\ge 0}.
\]

Let $\nu$ be a  2-form on $X$. We consider  formally the pairing
\begin{equation} \label{norma}
( f, g) \mapsto \int_{X}^{} f \bar{g} \nu + (f, g)_{1,X}.
\end{equation}

\noindent If $\nu$ is
positive and continuous, then the pairing \eqref{norma}  is a non degenerate, hermitian inner product on $\Luno (X)$, endowing it 
 with a Hilbert space structure. We denote by $\Vert \cdot \Vert_{\nu}$ the associated norm. As $X$ is compact, the topology on $\Luno (X)$ induced by this norm does not depend on the choice of $\nu$. The subspace $C^\infty (X)$ is dense in $\Luno(X)$ under this topology.

We can also make $\Luno$ into a sheaf. For an open set $U \subseteq X$, we denote by $\nu_U$ the restriction to $U$ of the given continuous volume form. We define $\Luno(U)$ as the space of distributions on $U$ which are $L^2$ with respect to $\nu_U$ and such that $\partial f, \bar{\partial}f$ are $L^2$. Again, $\Luno(U)$ is independent of the choice of $\nu$, as well as the choice of  ambient Riemann surface, cf \cite{Bost}, p.255. 

We define $\esp (X)$ to be the space of  $2$-currents $\mu$ that can be locally written on open sets $U$ as $\mu= dd^c f $ with $f \in \Luno(U)$. 

\begin{rem}
This definition  of $\esp$  is  equivalent to the one given in the introduction (cf. Theorem \ref{local} in the appendix).
\end{rem}

Let $\Sch(U)$ denote the space of $C^\infty$, complex valued functions on $X$ whose support is a compact set contained in $U$. We define $\Luno(U)_0$ to be the closure of the image of the natural map $\Sch (U) \rightarrow \Luno(X)$. As $X$ is compact, we have that $\Luno(X)_0 = \Luno(X)$. 

Similarly, we define $L^1(U)$ as the space of distributions on $U$ which are $L^1$ with respect to $\nu_U$, with $L^1(U)_0$ the closure of the image of the natural map $\Sch (U) \rightarrow L^1(X)$.

We  define formally $$(f, g)_{1,U} := i  \int_{U}^{} \partial f \wedge \overline{\partial g}.$$

\noindent Let us record for latter use the following result.

\begin{lema}\label{identidadutilisima} 
Assume either of the following hypotheses:

\begin{itemize}

\item $f \in \Luno(U)_0, g \in \Luno(U)$

\item $f \in L^1(U)_0, \frac{\partial }{\partial z} f \in L^1(U), g \in C^\infty(U)$.  

\end{itemize}

\noindent Then,  we have an equality of absolutely convergent integrals:

\begin{equation} \label{objetivo}
2\pi \int_{U}^{} f dd^c\bar{g} =- (f,g)_{1,U}.
\end{equation}

\end{lema}

\pruf It suffices to consider $f \in \Sch(U)$. Then, the identity $d(f\overline{\partial g}) = \partial f \wedge \overline{\partial g} - 2\pi i f dd^c \bar{g}$ and the vanishing of $f$ on $\partial U$ prove \eqref{objetivo}. Since either of the hypotheses ensure that $(f,g)_{1,U}$ is absolutely convergent, this finishes the proof. $\blacksquare$ \\

\subsection{Pull-back and push-forward in the space $\Luno$}

Let $\varphi : X \rightarrow Y$ be a non-constant holomorphic map between compact Riemann surfaces. For a function $g \in C^\infty (X )$ (resp. $f \in C^\infty (Y)$), the functions $\varphi^*f$ and $\varphi_* g$ are given by $\varphi^* f = f \circ \varphi$ and 

\begin{equation}\label{defpushforward}
\varphi_* g (y) = \sum_{x \in \varphi^{-1}(y)}  e_x g(x),
\end{equation}

\noindent where $e_x$ is the ramification index of $\ph$ at $x$. We have that $\varphi^* f \in C^\infty(X)$. Although  $\varphi_* g $ is continuous, it is not $C^\infty$ at branched points. However, we have the following result for the latter function.

\begin{lema}\label{tedioso}
Let $U\subset X$ be an open set.  For any $g \in C^\infty(U)$ the function $\varphi_* g $ given by \eqref{defpushforward} belongs to $ \Luno\big(\ph(U)\big)$.

\end{lema}

\pruf As $\varphi_* g $ is continuous, it suffices to check that the derivative of $\varphi_* g $ is $L^2$ around  the finite set of branched points contained in $\ph(U)$. Using local coordinates, we reduce the problem further to the following situation: $g$ is a $C^\infty$ function on $\disco = \{ z \in \ce, |z| <1 \}$, and  $\ph : \disco \rightarrow \disco$ is the map $\ph(z)=z^n$. Fixing a branch of the logarithm on $\disco$, we then have that $$\ph_* g(z) = \sum_{\rho^n=1}^{} g(\rho z^{1/n}) \quad \textrm{almost everywhere.}$$ Then, $$\frac{\partial \ph_*g}{\partial z}(z) = \frac{z^{(1/n)-1}}{n} h(z^{1/n}) \quad \textrm{ almost eeverywhere, where } h(z):= \sum_{\rho^n=1}^{} \rho  \frac{\partial g}{\partial z}(\rho z).$$ Hence,  

\begin{eqnarray}
 i\int_{\disco}^{}\Big| \frac{\partial \ph_*g}{\partial z}(z)\Big|^2\dz&=&\frac{i}{n^2}\int_{\disco}^{}|z|^{(2/n)-2} | h( z^{1/n})|^2\dz \nonumber \\
 &=& 2\int_{0}^{2\pi/n}\int_{0}^{1}r |h(re^{i\theta})|^2 dr d \theta,\nonumber
\end{eqnarray}

\noindent as can be seen by an elementary change of variables. The last integral converges because $h$ is smooth. $\blacksquare$ \\

\begin{prop} \label{continuidad}
The linear maps

\[
\varphi^* : C^\infty(Y) \rightarrow C^\infty(X), \quad \varphi_* : C^\infty(X ) \rightarrow \Luno(Y)
\]

\noindent are continuous in the $\Luno$ topology.

\end{prop}

We postpone the proof of this statement to the end of this Section. 

As the spaces $C^\infty(X),  C^\infty(Y)$ are dense in the respective $\Luno$ topologies,  Proposition \ref{continuidad} allows us to extend the pull-back and push-forward to the respective $\Luno$ spaces,

\[
\varphi^* : \Luno(Y) \rightarrow \Luno(X), \quad \varphi_* : \Luno(X ) \rightarrow \Luno(Y).
\]

\noindent We also have well-defined maps

$$\varphi^* : \esp(Y)  \rightarrow   \esp(X),   \quad \varphi_* :  \esp(X ) \rightarrow  \esp(Y),$$

\noindent  given locally by $\ph^* dd^c f := dd^c \ph^*f$ and $\ph_* dd^c g := dd^c \ph_*g$. \\

\begin{lema} \label{derivadafuntorial}

Fix an open set $U \subset X$. Assume either of the following hypotheses:

\begin{itemize}

\item $f \in \Luno\big(\ph(U)\big), g \in \Luno(U)$

\item $\frac{\partial}{\partial z}\ph^*f \in L^1(U), \frac{\partial}{\partial z} f \in L^1\big(\ph(U)\big), g \in C^\infty(U).$ 

\end{itemize}

\noindent Then,  we have that $$(\ph^*f,g)_{1,U} = (f,\ph_* g)_{1,\ph(U)}.$$

\end{lema}
\pruf Either of the hypotheses above allows us to reduce to the case where $f$ and $g$ are $C^\infty$  (by a well known density argument).  Then, using a partition of unity, the problem is further reduced to the case where $\ph$ is the covering $z \mapsto z^n$ of the complex unit disc. In this situation, the asserted equality  boils down to an elementary change of variables. $\blacksquare$\\

For a 2-current $\mu$ and a distribution $f$, we define  formally the symbol 

\begin{equation} \label{bilineal}
f\cdot \mu := \int_{X}^{}f \mu. 
\end{equation}

\begin{lema} \label{espfuntorial}

The integral in \eqref{bilineal} is absolutely convergent under either of the following hypotheses:

\begin{itemize}

\item[(i)] $f \in \Luno(X), \mu \in \esp(X)$ 

\item[(ii)] $f \in L^1(X), \frac{\partial}{\partial z} f \in L^1(X),  \mu $ a $C^\infty$ 2-form.

\end{itemize}

\noindent Moreover, we have that

\begin{enumerate}

\item \label{coninversa} $\ph^* s \cdot \mu = s \cdot \ph_* \mu$, for all $s \in \Luno(Y)$ and $\mu \in  \esp(X)$

\item \label{condirecta}$ \ph_* t \cdot \nu = t \cdot \ph^* \nu$, for all $t \in \Luno(X)$ and $\nu \in  \esp(Y)$.

\end{enumerate}
\end{lema}

\pruf Suppose hypothesis (i) holds. Using a partition of unity and the compactness of $X$, the integral in \eqref{bilineal} splits as a finite sum of integrals over open sets $U$ such that $\mu=dd^c h$, with $h \in \Luno (U)$ and   $f \in \Luno( U)_0$. Then, every integral on the sum converges by Lemma \ref{identidadutilisima}, which shows that \eqref{bilineal} is absolutely convergent. The same argument works under hypothesis (ii), taking $h \in C^\infty(U)$ and  $f \in L^1(U)_0$.

To prove \eqref{coninversa}, we choose a finite collection of open sets $U$ as in the paragraph above. As $\ph$ is a proper map,  we may suppose that $\ph^*s \in \Luno(U)_0$. Using  Lemma \ref{identidadutilisima} and Lemma \ref{derivadafuntorial}, we obtain that

\begin{eqnarray}
2\pi \int_{U}^{} \alpha^*s dd^c h & = & - (\alpha^*s,\bar{h})_{1,U} \nonumber \\
&=& - (s, \alpha_* \bar{h})_{1, \ph(U)} \nonumber \\
&=& 2 \pi \int_{\ph(U)}^{} s dd^c\alpha_* h.
\end{eqnarray}

The identity \eqref{condirecta} is proved similarly. $\blacksquare$\\

Now, we establish three lemmas leading to the proof of Proposition \ref{continuidad}. The first one will be used to control the $L^2$ norm of $\ph^*$ (resp. $\ph_*$) away from the ramification (resp. branched) points.

\begin{lema} \label{comparacion}
 Let $\mu_X$ (resp. $\mu_Y$) be a nonnegative, continuous 2-form on $X$ (resp. on $Y$). Let $\varphi : X \rightarrow Y$ be a non constant holomorphic map between compact Riemann surfaces, and let $f \in C^\infty (Y), g \in   C^\infty (X)$. Then, the following relations hold: 

\begin{eqnarray}
\| \varphi^* f \|_{\mu_X} & \leq & \sqrt{\deg \varphi} \| f \|_{\varphi_*\mu_X} \nonumber \\
\| \varphi_* g \|_{\mu_Y} & \leq &\deg \varphi \| g \|_{\varphi^*\mu_Y}. \nonumber 
\end{eqnarray}

\end{lema}

\begin{rem}

We stress that, as $\varphi_*\mu_X$ is \emph{not}  continuous if $\ph$ is ramified,  this lemma does not immediately imply Proposition \ref{continuidad} (cf. Remark \ref{usando DL}).

\end{rem}

\pruf  The space of $C^\infty$ functions is dense in $\Luno$. In particular, we can approximate in the $\Luno$ sense continuous functions by $C^\infty$ functions. Hence, we may suppose without
loss of generality that $\mu_X \in \esp(X)$ and $\mu_Y \in \esp(Y)$. We use Lemma 
\ref{espfuntorial} and Lemma \ref{derivadafuntorial} to find that

\begin{eqnarray}
\Vert \ph^*f \Vert_{\mu_X}^2 &=& |\ph^*f|^2 \cdot \mu_X  + (\ph^*f, \ph^*f)_{1,X}\nonumber \\
&=& \ph^* |f|^2 \cdot \mu_X + (f,\ph_* \ph^*f)_{1,Y} \nonumber \\
&=& |f|^2 \cdot \ph_* \mu_X + \deg \ph (f,f)_{1,Y} \nonumber \\
&\leq& \deg \ph \Vert f \Vert^2_{ \ph_* \mu_X}. \nonumber
\end{eqnarray}

This proves the first assertion. We will now deduce the second from the first:

\begin{eqnarray}
\Vert \varphi_* g \Vert^2_{\mu_Y} & = &  \varphi_* g  \cdot \big((\overline{\ph_* g})\mu_Y \big)+ (\ph_*g, \ph_*g)_{1,Y} \nonumber \\
& = &  \varphi_* g  \cdot \big( (\ph_* \overline{g}) \mu_Y\big) + (\ph_*g, \ph_*g)_{1,Y} \nonumber \\
&=&  g \cdot \ph^*  \big((\ph_* \overline{g} ) \mu_Y\big) + (g, \ph^* \ph_*g)_{1,X} \nonumber \\
&=&  g \cdot \big(\ph^* \ph_*\overline{g} \big) \big( \ph^*\mu_Y \big)+ (g, \ph^* \ph_*g)_{1,X}. \nonumber 
\end{eqnarray}

Applying  the Cauchy-Schwarz inequality to the inner product $$(h,l) \mapsto h\cdot \overline{l} \ph^* \mu_Y + (h,l)_{1,X},$$ we  obtain that

\begin{eqnarray}
\Vert \varphi_* g \Vert^2_{\mu_Y} & \leq & \Vert g \Vert_{\ph^*\mu_Y}  \Vert \ph^* \ph_*g \Vert_{\ph^*\mu_Y} \nonumber \\ 
&\leq & \sqrt{\deg \ph}  \Vert g \Vert_{\ph^*\mu_Y}  \Vert  \ph_*g \Vert_{\ph_* \ph^*\mu_Y} \nonumber \\ 
&= & \sqrt{\deg \ph}  \Vert g \Vert_{\ph^*\mu_Y}  \Vert  \ph_*g \Vert_{(\deg  \ph) \mu_Y} \nonumber \\ 
&\leq& \deg \ph  \Vert g \Vert_{\ph^*\mu_Y}  \Vert  \ph_*g \Vert_{\mu_Y}. \nonumber 
\end{eqnarray}

\noindent Lemma \ref{tedioso} ensures that  $  \Vert   \ph_* g\Vert_{ \mu_Y}$ is a real number. Hence, by the above inequality, we conclude that $\Vert \varphi_* g \Vert_{\mu_Y} \leq \deg \ph \Vert g \Vert_{\ph^* \mu_Y}. \blacksquare$ \\

\begin{lema} \label{formalidad}
Let $\disco:=\{|z| < 1\} \subset \ce$ and let $S=\{ f \in C^\infty (\overline{\disco}) : f=0$ on $\partial \disco$\}. Then, for all $f \in S$, we have that 
$$\int_{\disco}^{} \Big| \frac{\partial f}{\partial z}  \Big|^2 \dz =  \int_{\disco}^{} \Big| \frac{\partial f }{\partial \bar{z}}  \Big|^2 \dz.$$ 
\end{lema}

\pruf  Using Lemma \ref{identidadutilisima}, we have that

\begin{equation} \label{invarianza}
\int_{\disco}^{} \Big| \frac{\partial f}{\partial z}  \Big|^2 \dz = \int_{\disco}^{} \partial f \wedge \overline{\partial f} = 2\pi i\int_{\disco}^{} f dd^c \bar{f}.
\end{equation}

\noindent Similarly,

$$ \int_{\disco}^{} \Big| \frac{\partial f }{\partial \bar{z}}  \Big|^2 \dz = \int_{\disco}^{} \overline{\overline{\partial} f} \wedge \overline{\partial} f = \int_{\disco}^{} \partial \bar{f} \wedge \overline{\partial \bar{f}}= 2\pi i \int_{\disco}^{} \bar{f}dd^c f. $$

\noindent  The identity 
$d(f  \overline{\partial f} + \bar{f} \partial f) = 2\pi i (\bar{f}dd^cf-f dd^c\bar{f})$, along with the vanishing on the boundary allow us to deduce the result. $\blacksquare$ \\

The following lemma can be seen as a  weighted Poincar\'e inequality. It will be used in the proof of Proposition \ref{continuidad} to control the $L^2$ norm of $\ph^*$ around the ramification points.

\begin{lema} \label{ramificacion}
Let $ S$ be the set defined in Lemma \ref{formalidad}. 
 For every $\delta>0$ there exists $C_\delta > 0$ such that for all $f \in S$,

\begin{equation} \label{fcomplex}
i\int_{\disco}^{} \frac{|f(z)|^2}{|z|^{2-\delta}} \dz \leq C_\delta i \int_{\disco}^{} \Big| \frac{\partial f}{\partial z}\Big|^2 \dz.
\end{equation}

\end{lema}

\pruf Observe that the integral on the left hand side is convergent. Indeed,

$$i\int_{\disco}^{} \frac{|f(z)|^2}{|z|^{2-\delta}} \dz = 2 \int_{0}^{2\pi} \int_{0}^{1} \frac{|f|^2}{r^{1-\delta}} dr d\theta \leq \frac{4\pi}{\delta} \sup_\disco |f|^2.$$

Let $\disco_\eps = \{ \eps < |z| < 1 \}$. We have that

\begin{equation}\label{truco}
\frac{2}{\delta} d \Big( \frac{|z|^\delta}{\bar{z}} |f|^2d \bar{z} \Big) = \frac{|f|^2}{|z|^{2-\delta}} \dz + \frac{2}{\delta} \frac{|z|^\delta}{\bar{z}} \frac{\partial |f|^2}{\partial z} \dz \quad \textrm{ on } \disco_{\eps}.
\end{equation}

\noindent We remark that  $$\lim_{\eps \rightarrow 0}\int_{\partial \disco_\eps}^{}  \frac{|z|^\delta}{\bar{z}} |f|^2d \bar{z} = \lim_{\eps \rightarrow 0} -i\eps^\delta \int_{0}^{2\pi} |f(\eps e^{i\theta})|^2 d\theta = 0.$$  Then, using \eqref{truco},

\begin{eqnarray}
 i\int_{\disco}^{}\frac{|f|^2}{|z|^{2-\delta}} \dz &=&  \lim_{\eps\rightarrow 0} i \int_{\disco_\eps}^{}  \frac{|f|^2}{|z|^{2-\delta}} \dz \nonumber \\
&=& \lim_{\eps \rightarrow 0}  -\frac{2}{\delta}i\int_{\disco_\eps}^{} \frac{|z|^\delta}{\bar{z}} \frac{\partial |f|^2}{\partial z} \dz - i \frac{2}{\delta} \int_{\partial \disco_\eps}^{}  \frac{|z|^\delta}{\bar{z}} |f|^2d \bar{z} \nonumber \\  
&=& -\frac{2}{\delta}i\int_{\disco}^{} \frac{|z|^\delta}{\bar{z}} \frac{\partial |f|^2}{\partial z} \dz. \nonumber
\end{eqnarray}

But 

\begin{eqnarray}
 \Big| \int_{\disco}^{} \frac{|z|^\delta}{\bar{z}} \frac{\partial |f|^2}{\partial z} \dz \Big|&\leq& i \int_{\disco}^{} |z|^{\delta-1} \Big|\frac{\partial |f|^2}{\partial z}\Big| \dz \nonumber \\
 &\leq&  i \int_{\disco}^{} |z|^{\delta-1} \Big| f \frac{\partial f}{\partial z}\Big| \dz    +  i \int_{\disco}^{} |z|^{\delta-1} \Big| f \frac{\partial f}{\partial \bar{z}}\Big| \dz, \label{acato}
\end{eqnarray}

\noindent where we have used the relation $\frac{\partial }{\partial z} |f|^2= \bar{f} \frac{\partial }{\partial z}f + f \overline{\frac{\partial }{\partial \bar{z}} f}$. We bound the first integral using the Cauchy-Schwarz inequality. That is,

\begin{eqnarray}
 i \int_{\disco}^{} |z|^{\delta-1} \Big| f \frac{\partial f}{\partial z}\Big| \dz &\leq&  \Big( i\int_{\disco}^{} \frac{|f|^2 }{|z|^{2-2\delta}} \dz \Big)^{1/2} \Big( i\int_{\disco}^{} \Big| \frac{\partial f}{\partial z} \Big|^2 \dz \Big)^{1/2} \nonumber \\
 &\leq& \Big( i\int_{\disco}^{} \frac{|f|^2 }{|z|^{2-\delta}} \dz \Big)^{1/2} \Big( i\int_{\disco}^{} \Big| \frac{\partial f}{\partial z} \Big|^2 \dz \Big)^{1/2}. \nonumber
\end{eqnarray}

\noindent Using Lemma \ref{formalidad}, the same bound applies to the second integral in \eqref{acato}. Putting everything together, we obtain that

$$i\int_{\disco}^{}\frac{|f|^2}{|z|^{2-\delta}} \dz \leq  \frac{4}{\delta}   \Big( i\int_{\disco}^{} \frac{|f|^2 }{|z|^{2-\delta}} \dz \Big)^{1/2} \Big( i\int_{\disco}^{} \Big| \frac{\partial f}{\partial z} \Big|^2 \dz \Big)^{1/2}.$$ 

\noindent This implies \eqref{fcomplex} with $C_\delta = (4/\delta)^2. \blacksquare$ \\

\begin{rem}
Lemma \ref{ramificacion} can be viewed as a variant of the Hardy inequality (\cite{Inequalities}, Theorem 327). In both cases, the main idea is to express the weighted $L^2$  norm of $f$ as an integral involving the derivative of $f$. In our case, this is done through the equality \eqref{truco}. 
\end{rem}

\noindent \textbf{Proof} of Proposition \ref{continuidad} :  the endomorphism of the unit disk $\disco$ given by $z \mapsto z^n$ pulls back $\dz$ to $n^2 |z|^{2(n-1)}\dz$. Hence, the form $\varphi^*\mu_Y$ on Lemma \ref{comparacion} is continuous, non negative, with zeroes at the ramification points. This proves the assertion about the push-forward in Proposition \ref{continuidad}.

Concerning the pull-back, we must prove that we can choose a positive, continuous 2-form $\mu=\mu(\ph, \mu_X)$ on $Y$  as well as a constant $C = C(\ph, \mu_X) > 0$, such that

\begin{equation}\label{acotamiento}
\Vert \varphi^* f \Vert_{\mu_X} \leq C \Vert f \Vert_\mu
\end{equation}

\noindent holds for all $f \in C^\infty(Y)$. 

Lemma \ref{derivadafuntorial} implies the equality $(\ph^*f, \ph^*f)_{1,X} = \deg \ph (f,f)_{1,Y}$, so it suffices to  bound the $L^2$ norm of $\ph^*f$.

 Let $E \subset Y$ be the set of branched points of $\ph$. We select  a convenient neighborhood $U$ of $E$ as follows: for each ramification point  $P \in X $ of $\varphi$ with ramification index $n_P$, there exists a neighborhood $V_P \subset X$ (resp. $U_P \subset Y$) of  $P$ (resp. of the branched point $\varphi(P)$) and local charts $\theta_P :  \disco \rightarrow V_P $, $\eta_P : \disco \rightarrow U_P$ such that the following diagram commutes:

\[
\xymatrix{
 z \ar@{|->}[d] & \disco  \ar[d] \ar[r]^{\theta_P} &V_P  \ar[d]^{\varphi}    \\
z^{n_P} & \disco \ar[r]^{\eta_P} &  U_P 
}
\]

\noindent We may assume without loss of generality (shrinking the neighborhoods if required) that $\overline{V}_P \cap \overline{V}_Q = \emptyset$ if $P \neq Q$.
Let us put $U:= \cup_{P \textrm{ ramified}} U_P$. 

Let $U_E$ be another neighborhood of $E$ such that $\overline{U}_E \subset U$. Let $\psi_1, \psi_2$ be a partition of unity subordinate to the covering $\{ Y \backslash \overline{U}_E, U \}$. Let 

$$M:= \max \Big\{ \sup_{i=1,2} |\psi_i|,  \sup_{i=1,2} \Big| \frac{\partial}{\partial z} \psi_i \Big|, \sup_{i=1,2} \Big| \frac{\partial}{\partial \bar{z}} \psi_i \Big| \Big\}.$$

We  split   $f \in C^\infty(Y)$ as $f = f_1 + f_2$,  where  $f_1 := f \psi_1$ and $f_2:= f \psi_2$. As $$|f_i| \leq M |f|, \quad \Big| \frac{\partial  f_i}{\partial z} \Big| \leq  M \Big( |f| +   \Big|\frac{\partial f}{\partial z}\Big|\Big), \quad \Vert \ph^*f \Vert_{\mu_X} \leq \Vert \ph^*f_1 \Vert_{\mu_X} + \Vert \ph^*f_2 \Vert_{\mu_X},$$ it suffices  to prove \eqref{acotamiento} for $f_1$ and $f_2$.

Let $\mu$ be  a positive continuous volume form such that $\mu = \ph_*\mu_X $ on $Y \backslash  \overline{U}_E$. Note that such a form exists because $\ph_*\mu_X$ is continuous on 
$Y\backslash E$. We have that $\Vert f_1 \Vert_{\ph_*\mu_X} = \Vert f_1 \Vert_{\mu}$. Then, since supp$f_1 \subset Y \backslash \overline{U}_E$, it follows from Lemma \ref{comparacion} that 

$$\Vert \ph^*f_1\Vert_{\mu_X}\leq  \sqrt{\deg \ph} \Vert f_1 \Vert_{\mu}.$$

 We may suppose that the support of $f_2$ is contained in some $U_P$. We then have that

\begin{eqnarray}
 \int_{X}^{} |\varphi^*f_2 |^2 \mu_X & = &  \int_{\varphi^{-1}(U_P)}^{} |\varphi^*f_2 |^2  \mu_X \nonumber \\
& = &  \sum_{\ph(Q)=P}^{}\int_{\disco}^{} |\theta_Q^*\varphi^*f_2 |^2  \theta_Q^*\mu_X \nonumber \\
& \leq & \lambda \sum_{\ph(Q)=P}^{} i \int_{\disco}^{}|\theta_Q^*\varphi^*f_2 |^2 \dz. \nonumber 
\end{eqnarray}

\noindent The last inequality, valid for some $\lambda > 0$ depending only on the choice of charts, follows from the continuity of $\theta_Q^*\mu_X$ on $\overline{\disco}$. Putting $n:=n_P$ we then find that

\begin{eqnarray}
 i\int_{\disco}^{}|\theta_Q^*\varphi^*f_2 |^2 \dz &=& i \int_{\disco}^{} |f_2 \circ \eta_Q(z^{n})|^2 \dz \nonumber \\
 &=& \frac{i}{n}\int_{\disco}^{} \frac{|f_2 \circ \eta_Q (z)|^2  }{|z|^{2-(2/n)}} \dz \nonumber \\
&\leq& C_n i \int_{\disco}^{} \partial (f_2 \circ \eta_Q) \wedge \overline{\partial (f_2 \circ \eta_Q)}, \nonumber
\end{eqnarray}

\noindent where the last follows by application of Lemma \ref{ramificacion} to $f_2 \circ \eta_Q \in \Sch (\disco)$. Finally,

\begin{eqnarray}
\int_{\disco}^{} \partial (f_2 \circ \eta_Q) \wedge \overline{\partial (f_2 \circ \eta_Q)} &= &\deg \ph  \int_{U_P}^{} \partial f_2 \wedge \overline{\partial f_2} \nonumber \\
&=& \deg \ph \int_{Y}^{} \partial f_2 \wedge \overline{\partial f_2}. \nonumber 
\end{eqnarray}

This achieves the proof of \eqref{acotamiento} for $f_2$.  $\blacksquare$ \\

\begin{rem} \label{usando DL}
J.-B. Bost explained to us another way to handle the problem of the ramification points. Let $C_\ph^\infty(Y)$ be the subspace of $C^\infty(Y)$ consisting of functions whose support does not contain  the set of branched points of $\ph$. A result by Deny and Lions (\cite{DeLi}, 
Th\'eor\`eme II.2.2), using Cartan's theory of \emph{balayage}, ensures that   $C_\ph^\infty(Y)$ is dense in $\Luno(Y)$. Then to obtain a well defined pull-back map $\ph^* : \Luno(Y) \rightarrow \Luno(X)$, it suffices to  show the continuity of $\ph^*$ on $C_\ph^\infty(Y)$, which can be easily done (e.g. using Lemma \ref{comparacion} and the classical Poincar\'e inequality). This approach avoids the use of Lemma \ref{ramificacion}, though is less  elementary.
\end{rem}

\subsection{Pull-back and push-forward of $\Luno$-Green functions}

Let $X$ be a compact Riemann surface, and let $D = \sum_{P} a_P [P]$ be a divisor with real coefficients on $X$. By definition, a Green function with $C^\infty$ regularity for $D$ is a real-valued distribution $h$ on $X$ such that the current $dd^c h + \delta_D$ is a $C^\infty$ 2-form\footnote{This is what Gillet and Soul\'e call a \emph{Green current}. Following Bost, we use a different terminology better suited to the discussion of regularity issues. }. 

\begin{rem} \label{caracterizacion}
The above definition is equivalent to the following two properties: 

\begin{itemize}

\item $h$ is $C^\infty$ outside of the finite set $|D|:=\{ P $ such that $a_P \neq 0\}$

\item for every $P \in |D|$, there exists a local chart $(U, \theta_P)$ centered at $P$ such that

\[
 h = -a_P\log|\theta_P(\cdot)|^2 + b, 
\]

\noindent on $U$, where $b \in C^\infty(U)$.

\end{itemize}
The first property is a consequence of the ellipticity of the $dd^c$ operator. The second one comes from the fact that $\log|\cdot|^2$ is a fundamental solution of $dd^c f = \delta_0$ on the complex unit disc. 

\end{rem}

\begin{definicion}
An $L^2_1$-Green function for $D$ is a real-valued distribution $g$ on $X$ such that there exists a Green function with $C^\infty$ regularity $h$ for $D$, along with an element $\psi \in L^2_1(X)$ such that

\begin{equation} \label{defLuno}
g = h + \psi.
\end{equation}

\noindent We call this equality a \emph{regular splitting} of $g$.
\end{definicion}

\begin{rem}\label{medida}
A distribution $g$ satisfies the above definition if and only if $dd^cg + \delta_D$ is a 2-current in $\esp(X)$ (cf. Proposition \ref{existencia} below). In particular, the 2-current $dd^cg$ is a Radon measure, that is, it extends to a continuous functional on the space $C^0(X)$ of continuos functions (and not merely on $C^\infty(X)$).
\end{rem}

 Let $\varphi : X \rightarrow Y$ be a non-constant holomorphic map between compact Riemann surfaces. Let $D$ be a divisor on $X$ and let $h$ be a Green function 
with $C^\infty$ regularity for $D$. We define the distribution $\varphi_* h$ on $Y$ by
 
 \begin{equation} \label{additive}
\varphi_* h (y) := \sum_{x \in \varphi^{-1}(y)}^{} e_x h(x), 
\end{equation}

 \noindent where $e_x$ is the ramification index of $\varphi$ at $x$.

The definition above has to be understood in the following sense. Using Remark \ref{caracterizacion}, we see that the right side of (\ref{additive}) defines a function on the complement of the finite set

\[
\{ y \in Y \textrm{ such that } \varphi^{-1}(y) \cap |D| \neq  \emptyset \}. 
\] 

It can be proven that the singularities of $\varphi_* h$ on this set are logarithmic  (cf. \eqref{tipo logaritmico} below for a precise verification), thus defining a distribution on $X$. Finally, we define
 
  \[
\varphi_* g  :=  \varphi_*h + \varphi_*\psi.
\]

\noindent Similarly, let $E \in $Div$(Y)$ and let $k$ be a Green function with $C^\infty$ regularity for $E$. We define the distribution $\varphi^*k $ on $X$ by
 
 \[
\varphi^* k:= k \circ \varphi. 
\]

\noindent This definition is made rigorous in the same way as we did it with $\varphi_* h$ above.

\noindent  Let $f$ be a $\Luno$-Green function for $E$, with a regular splitting 
  
 \begin{equation} \label{defLuno2}
f = k + \xi,
\end{equation}

\noindent where $\xi \in \Luno (Y)$. We define 

\[
\varphi^* f := \varphi^*k + \varphi^* \xi.
\]

\begin{thm} \label{mainanalitico}
The distribution $\varphi_* g$ (resp. $\varphi^* f$) is a $\Luno$-Green function for $\varphi_* D$ (resp. for $\varphi^* E$).
\end{thm}

\pruf We first prove the assertion concerning the push-forward. It suffices to consider the case $D = [P]$. Let $n$ be the ramification index of $\varphi$ at $P$, and let $Q := \varphi(P)$. We choose local charts $(U, \theta_P)$ centered at $P$ and $(V, \theta_Q)$ centered at $Q$ such that the following diagram commutes:

\begin{equation} \label{discos}
\xymatrix{
U \ar[r]^{\theta_P} \ar[d]^{\varphi} & \disco \ar[d]^{\alpha} \\
V \ar[r]^{\theta_Q}  & \disco. }
\end{equation}

\noindent Here,  $\alpha(z) = z^n$.  We may assume that the ramification index is 1 for every point in $U$ different from $P$, and that we have

\[
h = -\log|\theta_P|^2 + b
\] 

\noindent on $U$, with $b \in C^\infty (U)$. We may also assume that supp $h \subset U$\footnote{If this is not the case, we take a $C^\infty$ function $s$ with support contained in $U$ and such that $s \equiv 1$ in a neighborhood of $P$. Then, $g= sh + \big( \psi + (1-s) h \big)$ is also a regular splitting as in \eqref{defLuno}.}. 

\noindent As $\ph_*\psi \in \Luno(Y)$, it suffices to prove that $\varphi_* h$ can be expanded around $Q$ as  

\[
\varphi_* h = - \log|\theta_Q(\cdot)|^2 + \tilde{b},
\]

\noindent where $\tilde{b} \in \Luno\big(\ph(U)\big)$. Using our assumptions on the support of $h$, we have that for all $y \in V  \backslash \{ Q\}$, 

\begin{eqnarray}
\varphi_* h (y)  & = &  \sum_{x \in \varphi^{-1}(y) \atop x \in U}^{} h(x) \nonumber\\
 & = &  \sum_{x \in \varphi^{-1}(y) \atop x \in U}^{} -\log|\theta_P(x)|^2 + b(x) \nonumber\\
&= &  \Big(\sum_{x \in \varphi^{-1}(y) \atop x \in U}^{} -\frac{1}{n}\log|\alpha \circ \theta_P(x)|^2\Big)  +\varphi_* b(y)  \nonumber \\
&= & \Big(\sum_{x \in \varphi^{-1}(y) \atop x \in U}^{} -\frac{1}{n}\log|\theta_Q \circ \varphi (x)|^2\Big)  +\varphi_* b(y) \nonumber \\
&= & \Big( \sum_{x \in \varphi^{-1}(y) \atop x \in U}^{} -\frac{1}{n}\log|\theta_Q (y)|^2 \Big)  +\varphi_* b(y)\nonumber \\
&= & -\log|\theta_Q (y)|^2 +\varphi_* b(y) \label{tipo logaritmico}
\end{eqnarray}

\noindent Using Lemma \ref{tedioso}, we have $\tilde{b}:=\varphi_* b \in \Luno\big(\ph(U)\big)$.

\noindent Now, we prove the assertion concerning the pull-back. As before, we may assume that $E=[Q]$. We have that

\[
\varphi^* [Q] = \sum_{\varphi(P) = Q}^{} e_P [P].
\]

\noindent We fix $P \in |\varphi^* [Q]|$ and choose local charts $(U, \theta_P)$ centered at $P$ and $(V, \theta_Q)$ centered at $Q$ such that the diagram (\ref{discos}) commutes, with $\alpha(z)=z^{e_P}$. Let us write $f = k + \xi$ as in (\ref{defLuno2}). This reduces the problem  to a study of the expansion of $\varphi^* k$ around $P$. We have that

\begin{eqnarray}
\varphi^* k & = & k \circ \varphi \nonumber \\
& = & -\log |\theta_Q \circ \varphi|^2 + b \circ \varphi \nonumber \\
& = & -\log |\theta_P^{e_P} |^2 + \varphi^* b  \nonumber \\
& = & -e_P\log |\theta_P |^2 + \varphi^* b. \nonumber 
\end{eqnarray}

\noindent As $\ph^*b$ is $C^\infty$, this concludes the proof. $\blacksquare$ \\

The following two Lemmas will be used in the proof of Theorem \ref{main}.

We adopt the following notation. For a  2-current $T$ and a test  function $\phi \in C^\infty$, we denote by $T[\phi]$ the value of $T$ at $\phi$.

\begin{lema} \label{preadmisibilidad} Using the notations of Theorem \ref{mainanalitico}, we have that

\begin{eqnarray} 
dd^c \ph^*f [ \phi ] &= & dd^c f [ \ph_* \phi ], \quad \phi \in C^\infty(X) \label{funtorialidadGreen} \\
dd^c \ph_*g[\phi ] & = & dd^c g [\ph^*\phi], \quad \phi \in C^\infty(Y). \label{laotrafuntorialidad}
\end{eqnarray}

\end{lema}

\begin{rem}
As $ \ph_* \phi$ is continuous,  the right hand side of \eqref{funtorialidadGreen} is well defined (cf. Remark \ref{medida}).

\end{rem}

\pruf we write $f= k + \xi$. Using Lemma \ref{espfuntorial}, we have that $$dd^c\ph^*\xi [\phi] =- \int_{X}^{} \ph^* \xi dd^c\phi = -\int_{Y}^{} \xi dd^c(\ph_*\phi) = dd^c\xi[\phi].$$

\noindent We put $w:= dd^c k + \delta_{E}$, which is a $C^\infty$ 2-form.  Using Lemma \ref{espfuntorial}, we have that

\begin{eqnarray}
dd^c \ph^*k [\phi] &=& \ph^*w[\phi] - \ph^*\delta_E[\phi] \nonumber \\
&=& \ph^*w \cdot \phi - \delta_E[\ph_*\phi] \nonumber \\
&=& w \cdot \ph_* \phi - \delta_E[\ph_*\phi] \nonumber \\
&=& dd^c k [\ph_*\phi], \nonumber 
\end{eqnarray} 

\noindent thus proving \eqref{funtorialidadGreen}. Equation $\eqref{laotrafuntorialidad}$ is proved similarly. $\blacksquare$\\

The following lemma is a variant of Lemma \ref{espfuntorial} for Green functions.

\begin{lema} \label{starfuntorial}
Let $h$ be a Green function with $C^\infty$ regularity for a divisor $D$ on $Y$ such that the support of $\ph^*h$ is contained in a neighborhood $U$ of $|\ph^*D|$. Let $s \in C^\infty(U)$. Then 
$$\int_{U}^{} \ph^* h dd^c s =\int_{\ph(U)}^{} h dd^c \ph_* s.$$
\end{lema}

\pruf Theorem \ref{mainanalitico} ensures that $\ph^*h$ and $\frac{\partial}{\partial z} \ph^*$ are $L^1$  on $U$ (i.e., apply  the local expansion in Remark \ref{caracterizacion} to $\ph^*h$). Then, by Lemma \ref{identidadutilisima} and Lemma \ref{derivadafuntorial}, we have that

\begin{eqnarray}
2\pi \int_{U}^{} \ph^* h dd^c s &=&  -(\ph^* h , \bar{s})_{1,U} \nonumber \\
 &=& -(h, \ph_*\bar{s})_{1,\ph(U)} \nonumber \\
 &=& 2\pi \int_{\ph(U)}^{}h dd^c \ph_* s. \blacksquare \nonumber
\end{eqnarray}

\section{Pull-back and push-forward in Arakelov theory}

\subsection{The $\Luno$-arithmetic Chow group, the admissible arithmetic Chow group and intersection theory} \label{Chowgroup}

Let $\X$ be an arithmetic surface, that is, a projective, integral, normal scheme over Spec$(\z)$ of Krull dimension 2. The ring $O_\X(\X)$ is the ring of integers of a number field $K$. We denote by Div$(\X)$ the set of Weil divisors on $\X$. 

A compactified divisor (in the sense of Bost) is a pair $(D, g)$, where $D \in $Div$(\X)$  and $g$ is a $L^2_1$-Green function for the divisor $D(\ce)$ on the compact Riemann surface $\X (\ce)$. In addition, $g$ must be invariant under complex conjugation on $\X (\ce)$. A principal divisor is a compactified divisor of the form $(\dv f, -\log |f_{\ce}|^2)$, where $f$ is a rational function on $\X$. This is a well defined notion, as can be seen from  the  Poincar\'e-Lelong equation (\cite{GriHa}, p. 388) 

\[
dd^c ( -\log |f_{\ce}|^2) + \delta_{\dv f_{\ce}}  = 0. 
\]

Let $\widehat{\textrm{Div}}(\X)$ be the set of compactified divisors. It is  an abelian group under the operation $(D,g) + (D',g') = (D+D', g+g')$. From remark  \ref{caracterizacion} 
it is easily verified that this operation is well defined. The neutral element is $(0, 0)$, that is, the empty divisor together with the zero function. The set of principal divisors is a subgroup of the group of compactified divisors.

The $\Luno$-arithmetic Chow group $\CH (\X)$ is defined as the quotient of the group of compactified divisors by the subgroup of principal divisors.

Let $X=\X(\ce)$. Let $\mu$ be a real-valued 2-current in $\esp(X)$ such that \begin{equation}\label{normal}
\int_{X}^{}\mu = 1. 
\end{equation} 

 A compactified divisor $(D, g)$ is said to be admissible with respect to $\mu$ (or $\mu$-admissible) if the following equality of currents holds:

\[
dd^c g + \delta_{D(\ce)} = \big(\deg D(\ce) \big) \mu.
\]

In this situation we will also say that $g$ is a $\mu$-admissible $\Luno$-Green function for $D$.

The $\Luno$-Green function $g$ is determined by $D$ up to a function which is constant on every connected component of $X$. Indeed, if $(D,h)$ is  $\mu$-admissible, then $g-h$ is an harmonic function on a compact Riemann surface. 
Hence, it must be constant on each connected component of $X$ by the maximum modulus principle.

We remark that  (\ref{normal}) is a necessary condition for the existence of admissible compactified divisors $(D,g)$ with $\deg D(\ce)\neq 0$. To give a  sufficient condition we use an existence result coming from the theory of elliptic operators: 

\begin{lema} \label{Poisson}
Let $\mu \in \esp(X)$. The equation 
\begin{equation} \label{ecuaciondePoisson}
dd^c u = \mu
\end{equation} has a solution $u \in \Luno(X)$ if and only if $\int_{X}^{}\mu =0$. If $\mu$ is real-valued then $u$ can be taken to be real-valued.
\end{lema}

If we take $\mu$ to be a $C^\infty$ form (in which case $u \in C^\infty(X)$), this is a well known consequence of the Hodge theorem  in differential geometry (cf. \cite{Warner}, 6.8).  

Equation \eqref{ecuaciondePoisson} can be  solved locally via standard techniques from the theory of elliptic PDE's, namely the Lax-Milgram theorem and the Fredholm alternative. To prove that it is possible to patch together these local solutions to get a global solution, we shall use the cohomology machinery. We explain the precise steps giving a proof of Lemma \ref{Poisson} in the Appendix.

\begin{prop}\label{existencia}
Let $\mu$ be a real-valued 2-current in $\esp(X)$  such that $\int_{X}^{}\mu =1$. Then for every Weil divisor $D$ on $\X$, there exists a $\Luno$-Green function $g$ such that $(D,g)$ is $\mu$-admissible.
\end{prop} 

\pruf let $\mu_{0}$ be the canonical Arakelov 2-form if the genus of $X$ is positive and the Fubini-Study form if it is zero.  Let $\tilde{g}$ be a real-valued solution to the equation of currents (\cite{Lang}, Theorem 1.4) $$dd^c \tilde{g} + \delta_{D(\ce)} = \big(\deg D(\ce)\big)\mu_{0}.$$ 

Let $\tilde{\mu} = \big(\deg D(\ce)\big)(\mu_{0} - \mu)$. Then $\tilde{\mu}$ belongs to $\esp(X)$ and $\int_{X}^{}\tilde{\mu} =0$. Using  Lemma \ref{Poisson}, we have  a real-valued $u\in \Luno(X)$ such that $dd^c u=\tilde{\mu}$. Then $g:=\tilde{g}-u$ is a $\mu$-admissible, $\Luno$-Green function for $D(\ce). \blacksquare$\\

Principal divisors  are admissible with respect to any 2-current. 
The set of $\mu$-admissible compactified divisors forms a group. We denote by $\CH(\X)_\mu$ the quotient of this group by the subgroup of principal divisors. The natural map $\CH(\X)_\mu \rightarrow \CH(\X)$ is an embedding of abelian groups.

Let $x_1,x_2 \in \CH(\X)$ be represented by the compactified divisors $(D_1,g_1)$ and $ (D_2, g_2)$. We suppose that the representatives have been chosen without  common irreducible components. The intersection pairing is defined as $$\izq x,y\der:=  \izq D_1,D_2 \der +\izq g_1,g_2\der,$$ where the first pairing is the finite intersection pairing as defined in \cite{Bost}, Section 5.3. We recall from loc. cit. the definition of the second pairing: we choose regular splittings $$g_i = h_i + \psi_i$$ as in \eqref{defLuno}, and we define
$w_i = dd^c h_i +
\delta_{D_i(\ce)}$. Then,

\begin{equation} \label{infinito}
\mathopen<g_1,g_2\mathclose> := \frac{1}{2} \bigg(\int_{X} h_1 * h_2 + \int_{X}
\psi_1 w_2 + \int_{X} \psi_2 w_1 + \frac{1}{2\pi i} \int_{X}
\partial \psi_1 \wedge \bar{\partial} \psi_2 \bigg),
\end{equation}

\noindent where $h_1 * h_2 = h_1w_2 + h_2\delta_{D_1(\ce)}$. This number does not depend of the choice of $h_i, \psi_i$.

\begin{ex} \label{grado cero}
For every compactified divisor $\hat{D}=(D,g)$ and $c \in \erre$, we have that $$\izq (0,c) , \hat{D} \der =\frac{c}{2}[K:\q]\deg D(K).$$
\end{ex}

Let $\ph : \X \rightarrow \Y$ be a generically finite morphism between arithmetic surfaces. If we have a regular splitting $g= h + \psi$, then by definition $\ph_* g = \ph_* h + \ph_*\psi$. As $\ph_* h$ is not a Green function with $C^\infty$ regularity when $\ph$ is ramified, this is not a regular splitting of $\ph_* g$. However, we show in the following lemma that this decomposition can still be used to compute intersections.

\begin{lema} \label{descomposicionauxiliar}
Let $(D_x,g_x)$ (resp. $(D_y,g_y)$) be a compactified divisor on $\X$ (resp.$\Y$). We decompose $g_x = h_x + \psi_x$ (resp. $g_y=h_y+\psi_y$) as in \eqref{defLuno} and we put $w(\ph_*h_x):=dd^c \ph_* h_x +\delta_{\ph_*D_x(\ce)}$ (resp. $w(h_y):=dd^ch_y + \delta_{D_y(\ce)}$). we have that 
$$\izq \ph_*g_x, g_y \der = \frac{1}{2} \bigg(\int_{Y} \ph_*h_x * h_y + \int_{Y}
\ph_* \psi_x w(h_y) + \int_{Y} \psi_y w(\ph_* h_x) + \frac{1}{2\pi i} \int_{Y}
\partial (\ph_*\psi_x) \wedge \bar{\partial} \psi_y \bigg).$$ Moreover, all of the integrals above are absolutely convergent. 
\end{lema}

\pruf Let $Y=\Y(\ce)$. Using Theorem \ref{mainanalitico}, we obtain an auxiliary regular splitting 
\begin{equation}\label{auxiliar}
\ph_* g_x = h + \psi,
\end{equation}

\noindent where $h$ is a Green function  for $\ph_*D_x$ with $C^\infty$ regularity, and  $\psi \in \Luno(Y)$. Hence, we have that

 $$\izq \ph_*g_x, g_y \der =  \frac{1}{2} \bigg(\int_{Y} h * h_y + \int_{Y}
\psi w(h_y) + \int_{Y} \psi_y w( h) + \frac{1}{2\pi i} \int_{Y}
\partial \psi \wedge \bar{\partial} \psi_y \bigg).$$ 

We deduce from \eqref{auxiliar} that

\begin{equation}\label{primera}
\int_{Y}^{} \ph_* h_x * h_y = \int_{Y}^{}h * h_y + \int_{Y}^{} (\psi - \ph_*\psi_x) w(h_y).
\end{equation}

 As $w(h_y)$ is $C^\infty$ and $\psi - \ph_*\psi_x$ is $L^2$, this equality shows the absolute convergence of the integral on the left. 
Thus, we find that
\begin{eqnarray}
\int_{Y} \psi_y w(\ph_* h_x)&=& \int_{Y} \psi_y w(h + \psi - \ph_*\psi_x) \label{paja} \\
 &=& \int_{Y}^{}\psi_y w(h) + \int_{Y}^{}\psi_y dd^c (\psi - \ph_* \psi_x)  \nonumber
\end{eqnarray}
 
The last integral is absolutely convergent because 

$$ \int_{Y}^{}\psi_y dd^c (\psi - \ph_* \psi_x) = \frac{1}{2\pi i}\int_{Y}^{} \partial \psi_y \wedge \overline{\partial}  (\psi - \ph_* \psi_x),$$ and the functions on the right are $\Luno$.
This shows that the left hand side of \eqref{paja} is absolutely convergent. From this equality we deduce that

\begin{equation}\nonumber
\int_{Y} \psi_y w(\ph_* h_x) + \frac{1}{2\pi i} \int_{Y}
\partial (\ph_*\psi_x) \wedge \bar{\partial} \psi_y  =   \int_{Y} \psi_y w( h) +  \frac{1}{2\pi i} \int_{Y}
\partial \psi \wedge \bar{\partial} \psi_y. 
\end{equation}

\noindent Together with \eqref{primera}, this proves both the absolute convergence of the first integral on the left and the desired formula. $\blacksquare$\\

\subsection{Functoriality with respect to pull-back and push-forward}

Let $\varphi: \X \rightarrow \Y$ be a generically  finite morphism between arithmetic surfaces.  Let $(D,g_D)$ (resp. $(E, g_E)$)  be a compactified divisor on $\X$ (resp. on $\Y$). We define 

\begin{eqnarray}
\varphi^* (E, g_E) & :=  & (\varphi^*E, \varphi^* g_E) \nonumber \\
\varphi_* (D, g_D) & :=  &(\varphi_*D, \varphi_* g_D). \nonumber 
\end{eqnarray}  

\begin{thm} \label{main}

\begin{enumerate}

\item \label{bien definido} The formulas above induce homomorphisms 
\begin{eqnarray}
\varphi^*  :   \CH(\Y) \longrightarrow \CH(\X)  \nonumber \\
\varphi_*  :   \CH(\X) \longrightarrow \CH(\Y).  \nonumber 
\end{eqnarray}

\item \label{funtorial} Let $\mu \in \esp \big(\X(\ce)\big)$ and $\nu \in \esp \big(\Y(\ce)\big)$ be real-valued and such that $\int_{\X(\ce)}^{}\mu = \int_{\Y(\ce)}^{}\nu=1$.  The formulas above induce homomorphisms

\begin{eqnarray}
\varphi^*  :   \CH(\Y)_\nu \longrightarrow \CH(\X)_{\frac{\varphi^*\nu}{\deg \ph}}  \nonumber \\
\varphi_*  :   \CH(\X)_\mu \longrightarrow \CH(\Y)_{\varphi_*\mu}.  \nonumber 
\end{eqnarray}  

\item \label{formula de proyeccion} Projection formula: for $x \in \CH (\X)$ and $y \in \CH(\Y)$, we have that
\[
\izq\varphi^* y, x\der = \izq  y, \varphi_*x\der.
\]
\end{enumerate}

\end{thm}

\begin{rem}
The functorial properties of the pairing in (\ref{bilineal}) ensure that $$\frac{1}{\deg \ph} \int_{X}^{} \ph^* \nu = 1 \textrm{ and }\int_{Y}^{} \ph_* \mu =1.$$ Hence, the arithmetic Chow groups appearing to the right in Theorem \ref{main} (\ref{funtorial}) contain compactified divisors of nonzero generic degree (cf. Proposition \ref{existencia}). 
\end{rem}
A correspondence $T$ on $\X$ is, by definition, an arithmetic surface $\Y$ together with two ordered finite morphisms $p,q : \Y \rightarrow \X$. The correspondence
 $T$ is said to be symmetric if $p_*q^* = q_* p^*$ on $Div(\X)$.

For a compactified divisor $(D, g)$ on $\X$ and $\mu \in \esp\big(\X(\ce)\big)$, we define

\begin{eqnarray}
\hat{T} (D, g) & := & (p_* q^* D, p_* q^* g) \label{Taritmetico} \\
T \mu & := & p_* q^* \mu \nonumber
\end{eqnarray}

The following statement is a straightforward consequence of Theorem \ref{main} and the definitions above. 

\begin{cor} \label{main en correspondencias}

\begin{enumerate}

Let $T$ be a  correspondence on the arithmetic surface $\X$. 

\item The formula (\ref{Taritmetico}) induces a homomorphism

\[
\hat{T} : \CH (\X) \longrightarrow \CH(\X).
\]

\item \label{simetrico} If $T$ is symmetric, then the morphism $\hat{T}$ is self-adjoint, i.e.  for all $x,y \in \CH(\X)$,

\[
\izq \hat{T}x, y \der = \izq x, \hat{T} y \der.
\]

\item Formula (\ref{Taritmetico}) induces a homomorphism

\[
\hat{T} : \CH (\X)_\mu \longrightarrow \CH(\X)_{\frac{T\mu}{\deg q}}.
\]

\end{enumerate}

\end{cor}

\begin{cor} \label{commutacion aritmetica}
 Let $S$ be another correspondence on $\X$.  If  $S$ and $T$ commute as endomorphisms of Div$(\X)$, then $\hat{S}$ and $\hat{T}$  commute as endomorphisms of $\widehat{\textrm{Div}}(\X)$.
\end{cor}

\pruf let $(D,g)$ be a compactified divisor. Let $g = h + \psi$ be a decomposition  as in \eqref{defLuno}. Let $\psi_n$ be a sequence of $C^\infty$ functions converging to $\psi$ in $\Luno\big(\X(\ce)\big)$. The hypothesis  implies that $S$ and $T$ commute as correspondences on $\X(\ce)$, so $ST\psi_n - TS\psi_n \equiv 0$. By Proposition \ref{continuidad}, we have that $ST \psi = TS \psi$.

By Theorem \ref{mainanalitico}, we deduce that
$$ST g - TS g = ST h - TS h \in \Luno\big(\X(\ce)\big).$$

The function $h$ is $C^\infty$ outside of the finite set $F:=\{ x \textrm{ such that } |ST [x]| \cap |D(\ce)| \neq \emptyset \} $. For any point $x \in \X(\ce) - F$, we have $ST h (x) = TS h (x)$ by hypothesis, so $ST h = TS h$ in $\Luno\big( \X(\ce)\big).$ This shows that $ST(D,g) = TS(D,g)$, as required. $\blacksquare$ \\

\noindent \textbf{Proof} of Theorem \ref{main}: To prove (\ref{bien definido}), it suffices by Theorem \ref{mainanalitico}  to  show that principal divisors are sent into principal divisors. Let $f$ be a rational function on $\X$. Let $N_\ph : K(\X) \rightarrow K(\Y)$ be the norm map between function fields induced by $\ph$.  We have that $$\ph_*(\dv f, -\log|f|^2)=(\dv N_\ph(f), -\log| N_\ph(f)|^2).$$ 

On the other hand, if $f$ is a rational function on $\Y$, we have that $$\ph^*( \dv f, -\log|f|^2) = (\dv \ph^*f, -\log| \ph^*f|^2).$$ This finishes the proof of (\ref{bien definido}).

Now we  prove the first statement of (\ref{funtorial}). Let $(E,g)$ be a compactified divisor on $\Y$ which is admissible with respect to $\nu$. Theorem \ref{mainanalitico} shows that $\ph^*(E,g)$ is a compactified divisor on $\X$. We  need to check the equality of currents $$dd^c (\ph^* g) + \delta_{\ph^*E(\ce)} = \deg \big(\ph^* E(\ce)\big)\frac{\ph^*\nu}{\deg \ph}.$$

 Let $\phi \in C^\infty(X)$.  Using Lemma \ref{preadmisibilidad}, we have that

\begin{eqnarray}
dd^c (\ph^* g)  [\phi]&= & dd^c g [\ph_*\phi] \nonumber \\
&=& -\delta_{E(\ce)} [\ph_* \phi]  + \big(\deg E(\ce)\big) \nu [\ph_* \phi] \nonumber \\
&=& -\delta_{\ph^* E(\ce)} [\phi]  + \big(\deg E(\ce)\big) \ph_* \phi\cdot \nu. \nonumber 
\end{eqnarray}

Using Lemma \ref{espfuntorial} and the fact that $\deg \ph^*E(\ce) = \deg \ph \deg E(\ce)$, this finishes the proof of the statement concerning $\ph^*$ in \eqref{funtorial}.

To prove the part concerning $\ph_*$ in \eqref{funtorial}, let $(D,g)$ be a compactified divisor on $\X$ which is admissible with respect to $\mu$. Theorem \ref{mainanalitico} shows that $\ph_*(D,g)$ is a compactified divisor on $\Y$. As $\deg \big(\ph_* D(\ce)\big) = \deg \big( D(\ce)\big) $, we  need to check the equality of currents $$dd^c (\ph_* g) + \delta_{\ph_*D(\ce)} = \deg \big( D(\ce)\big)\mu.$$ Let $\phi$ be a $C^\infty$ function on $Y$.  Using Lemma \ref{preadmisibilidad}, we have that

\begin{eqnarray}
dd^c (\ph_* g)[\phi] &=& dd^c g [\ph^* \phi] \nonumber \\
&=& -\delta_{D(\ce)} [\ph^* \phi] + \big(\deg D(\ce)\big) \mu[\ph^*\phi] \nonumber  \\
&=& -\delta_{\ph_*D(\ce)} [\phi] +  \big(\deg D(\ce)\big) \ph^*\phi \cdot  \mu. \nonumber  
\end{eqnarray}

As before, we conclude   by Lemma \ref{espfuntorial}. This finishes the proof of (\ref{funtorial}).

Now, we  prove (\ref{formula de proyeccion}). Let $(D_x,g_x)$ (resp. $(D_y, g_y)$) be a member in the class  $x$ (resp. $y$), chosen such that $|D_x(\ce)|\cap |\ph^*D_y(\ce)|=\emptyset$. We decompose $$g_x = h_x + \psi_x, \quad g_y=h_y + \psi_y,$$  as in (\ref{defLuno}). As $w(h_x)$ is a $C^\infty$ 2-form, there exists a neighborhood $U$ of $|\ph^*D_y(\ce)|$ and  $s \in C^\infty(U)$  with $w(h_x)=dd^c s$ on $U$. We suppose that the support of $\ph_*h_y$ is contained in $U$.   Using Lemma \ref{descomposicionauxiliar}, we have that $$\izq \ph^*g_y, g_x\der = \int_{X}^{} \ph^* h_y * h_x + \int_{X}^{} \ph^* \psi_y w(h_x) + \int_{X}^{} \psi_x w(\ph^* h_y) + \frac{1}{2\pi i} \int_{X}^{} \partial (\ph^* \psi_y) \wedge \bar{\partial} \psi_x.$$ We will prove the required equality term by term.  
 
 \begin{eqnarray} 
 \int_{X}^{} \ph^* h_y * h_x & = &\int_{X}^{} \ph^*h_y w( h_x)  + h_x\big(\ph^*D_y(\ce)\big) \nonumber\\
 &=& \int_{U}^{} \ph^*h_y w( h_x) +  \ph_*h_x\big(D_y(\ce)\big)  \nonumber\\
 &=& \int_{U}^{} \ph^*h_y dd^cs + \ph_*h_x\big(D_y(\ce)\big)  \nonumber\\
&=& \int_{\ph(U)}^{} h_y dd^c\ph_* s + \ph_*h_x\big(D_y(\ce)\big)  \label{acatonomas}\\
&=&  \int_{\ph(U)}^{} h_y \ph_* w(h_x) +  \ph_*h_x\big(D_y(\ce)\big)  \nonumber\\
&=& \int_{Y}^{} h_y *  \ph_* h_x, \nonumber
\end{eqnarray}

\noindent where we have used Lemma \ref{starfuntorial} in \eqref{acatonomas}. Using Lemma \ref{espfuntorial}, we obtain that
\begin{eqnarray}
 \int_{X}^{} \ph^* \psi_y w(h_x) &=& \int_{Y}^{}\psi_y \ph_* \big( w(h_x)\big) \nonumber \\
 &=&\int_{Y}^{}\psi_y  w(\ph_*h_x). \nonumber 
 \end{eqnarray}

The equality $ \int_{X}^{} \psi_x w(\ph^* h_y) = \int_{Y}^{} \ph_* \psi_x w (h_y)$ is handled similarly. Using Lemma \ref{derivadafuntorial} we obtain that

$$\int_{X}^{} \partial (\ph^* \psi_y) \wedge \bar{\partial} \psi_x= \int_{Y}^{}  \partial ( \psi_y) \wedge \bar{\partial} \ph_* \psi_x. \blacksquare$$\\

\section{Hecke correspondences on modular curves}

\subsection{Integral models of the modular curves $\XC{N}$ and Hecke correspondences}
Let $N \geq 1$ be an integer.  
We consider the modular curve $\XC{N}:= \G (N) \backslash \big(\h \cup \p(\q)\big)$. This Riemann surface admits a smooth, projective model over $\q$. In order to obtain a model over $\z$, we will use the modular interpretation in terms of $\G(N)$-structures à la Drinfeld as  in \cite{KatzMazur}. We recall the basic facts and definitions:

\begin{itemize}

\item Let $S$ be a scheme and let $E \rightarrow S$ be an elliptic curve. A group subscheme $G \subset E$ over $S$ is \emph{cyclic of order $N$} if it is locally free of rank $N$ and there exists a morphism $T \rightarrow S$ which is faithfully flat and locally of finite presentation and a $T$-point $P \in G (T)$ such that we have the equality of Cartier divisors $$G\times_S T = \sum_{a=1}^{N}[aP].$$

\item An isogeny  $\pi: E \rightarrow E'$ over $S$ is said to be cyclic of order $N$ if $\ker \pi$ is a cyclic subgroup of order $N$ in the above sense. 

\item Let $d$ be a divisor of $N$. There is a unique cyclic subgroup $G_d \subset G$, called the \emph{standard cyclic subgroup of order $d$}, characterized by:  if $P$ is a f.p.p.f.\footnote{Faithfully flat of finite presentation (\emph{Fidèlement plat, de présentation finie} in french). } local generator of $G$, then $G_d$ is generated by $(N/d)P$ (Theorem 6.7.2. of \emph{loc. cit.}).

\end{itemize}

Let $\XZ{N}$ be the compactified coarse moduli scheme associated to the moduli problem $[\G(N)]$ classifying cyclic $N$-isogenies between elliptic curves as constructed in \emph{loc. cit.} Ch. 8. Alternatively, the functor $[\G(N)]$ can be described as classifying cyclic subgroups of order $N$ of elliptic curves. The scheme $\XZ{N}$ is proper and normal. It has good reduction at $p \nmid N$ and bad reduction at $p |N$.

We consider the related functor $$F_N : \textrm{Sch} \longrightarrow \textrm{Sets}$$ which attaches to a scheme $S$ the set $$F_N(S) = \{ (E, C) \textrm{ such that } E/S \textrm{ is an elliptic curve and } C \textrm{ is a } \G(N)-\textrm{structure on }E \} / \sim,$$ where $\sim$ stands for the natural isomorphism notion. 

Let $l$ be a prime number such that $l \nmid N$. We define morphisms $$\alpha_l, \beta_l :  F_{Nl} \rightarrow F_N$$ as follows. Let $S$ be a scheme, and let $(E,C) $ represent a class in $F_{Nl}(S)$. 
Let $C_N \subset C$ (resp.  $C_l \subset C$) be the standard cyclic subgroup of order $N$ (resp. $l$).  We put

$$\alpha_l(E,C) := (E, C_N), \quad  \beta_l (E,C) := (E/C_l, C/C_l).$$ The morphism $\beta_l$ is well defined because of Theorem 6.7.4. of \emph{loc. cit.}

Using the coarse moduli property and the construction of the compactification as a normalization, we obtain finite morphisms $\alpha_l, \beta_l : \XZ{Nl} \rightarrow \XZ{N}$. We define the Hecke correspondence of order $l$ by $$T_l:= (\beta_l)_* \alpha_l^*.$$

Let $d$ be a divisor of $N$ such that $(d, N/d) = 1$. We define the Atkin-Lehner involution $w_d$ on $\XZ{N}$ first as a functor $F_N \rightarrow F_N$ as follows: let $S$ be a scheme and let $(E,C)$ represent a class in $F_N(S)$. Let  $C_d$ be the standard cyclic subgroup of order $d$   of $C$. We put  

\[
w_d (E, C) := \big(E/C_d, (E[N]\cap \frac{1}{d}C)/C_d\big).
\] 

As before, we obtain an involution $w_d : \XZ{N} \rightarrow  \XZ{N}$. \\

In what follows we will suppose $N$ squarefree to simplify the analysis of the bad fibers. For  a prime number $p$, we put $X_p(M):= \XZ{M} \otimes \Fp$. If $p|N$,  then $X_p(N)$ has two irreducible components, each one of them isomorphic to $X_p(N/p)$. We denote by $X_p^\infty(N)$ (resp. $X_p^0(N)$) the irreducible component intersecting $D_\infty^N$ (resp. $D_0^N$). 

 \begin{lema}\label{fibras finitas}
Suppose that $N$ is squarefree. For prime numbers $p|N$, $l \nmid N$ and  $u \in \{0, \infty\}$, we have that $$\alpha_l^*X_p^u(N) = \beta_l^*X_p^u(N)= X_p^u(lN),$$ 
$$ (\alpha_l)_*X_p^u(lN) = (\beta_l)_*X_p^u(lN) = (l+1)X_p^u(N),$$
$$ w_N \big(X_p^\infty(N) \big) = X_p^0(N), \quad w_N \big(X_p^0(N) \big) = X_p^\infty(N).$$
\end{lema}

\pruf  We first remark that $(\alpha_l)_* D_u^{Nl} = (\beta_l)_*  D_u^{Nl} = D_u^N$, and $w_N D_\infty^N = D_0^N$. (It is enough to check these equalities on the generic fiber, where it is clear). Using the notations $1/0 = \infty$ and $1/\infty =0$, we have that

\begin{eqnarray}
\izq w_N X_p^\infty(N), D_u^{N} \der &=& \izq X_p^\infty(N), w_N D_u^{N} \der \nonumber \\
 &=& \izq X_p^\infty(N),  D_{1/u}^{N} \der.  \nonumber
\end{eqnarray}

If $u=\infty$ (resp. $u=0$), then the last number is zero (resp. nonzero) by definition. This proves the assertion concerning $w_N$. Now, we have that 

\begin{eqnarray}
\izq \alpha_l^* X_p^\infty(N), D_u^{Nl} \der &=& \izq X_p^\infty(N), (\alpha_l)_* D_u^{Nl} \der \nonumber \\
 &=& \izq X_p^\infty(N),  D_u^{N} \der  \nonumber
\end{eqnarray}

If $u=\infty$ (resp. $u=0$), then the last number is nonzero (resp. zero) by definition, which shows that 

\begin{equation} \label{parte previa}
\alpha_l^*X_p^\infty(N) = X_p^\infty(Nl).
\end{equation}

 A similar argument applies to prove $\alpha_l^*X_p^0(N) = X_p^0(Nl)$. The analogous claims concerning $\beta_l^*$ are proved in the same way.

Now, for the push-forward, we have that

$$\izq (\alpha_l)_* X_p^\infty(Nl) , D_\infty^{N} \der = \izq  X_p^\infty(Nl) , \alpha_l^*D_\infty^{N} \der.  $$

Because $\alpha_l$ is finite, the divisor $ \alpha_l^*D_\infty^{N}$ is effective and horizontal. As  $D_\infty^{Nl}$ belongs to $\alpha_l^{-1}(D_\infty^N)$, the last number is nonzero. Using the fact that $ (\alpha_l)_* X_p^\infty(Nl)$ must  be irreducible, we conclude that it is supported on $X_p^\infty(N)$. 
To compute the degree of the extension of function fields over $\Fp$ given by $\alpha_l^* : K\big(X_p^\infty(N)\big) \rightarrow K\big(X_p^\infty(Nl)\big)$, we observe that the coarse moduli property together with \eqref{parte previa} allow to check that for every geometric point of $X_p^\infty(N)$, not supersingular, there are exactly $l+1$ points above it. Hence, the degree of $\alpha_l$ is $l+1$ and we have that $(\alpha_l)_* X_p^\infty(Nl) = (l+1) X_p^\infty(N)$. The same argument works with $0$ instead of $\infty$ or with $\beta_l$ instead of $\alpha_l. \blacksquare$

\begin{lema} \label{Heckesimetrico}
The correspondence $T_l$ is symmetric.
\end{lema}

\pruf this is well known on the generic fiber. Hence, for any irreducible, horizontal divisor $D$ on $\XZ{N}$, we have that $$(\beta_l)_*\alpha_l^* D = (\alpha_l)_*\beta_l^* D + V,$$ where $V$ is a vertical divisor. (This is justified because both sides of this equality are equal on the generic fiber). As $\alpha_l$ and $\beta_l$ are finite, no vertical divisor is produced during the pull-back operation. Hence,  $V=0$. On the other hand, if $D$ is an irreducible component of a finite fiber, we conclude by Lemma \ref{fibras finitas}. $\blacksquare$ \\

The cusps $\G(N)\infty$ and $\G(N)0$ induce $\q$-points of $\XZ{N}$. We denote by $D_\infty^N$ (resp. $D_0^N$) the corresponding Zariski closure in $\XZ{N}$.

\begin{lema} \label{dinfinito}
we have that $T_l D_\infty^N = (l+1)D_\infty^N$.
\end{lema}

\pruf An explicit calculation shows that the equality holds on the generic fiber. Then $T_l D_\infty^N - (l+1)D_\infty^N$ must be  a vertical divisor $V$. Since the morphisms $\alpha_l, \beta_l$ defining $T_l$ are finite, no vertical component is produced during the pull-back operation. Hence, $V=0. \blacksquare$ \\

\subsection{Arithmetic Chow group of $\XZ{N}$}

Let $\mu = \frac{dz \wedge d \bar{z}}{\big(Im (z)\big)^2}$ be the hyperbolic 2-form on $\h$.  This form is $GL_2^+(\erre)$-invariant, and  induces a finite measure on $\XC{N}$. We denote by $\mu_N:= \mu / Vol_{\mu}\big(\XC{N}\big)$ the normalized 2-form on $\XC{N}$.

\begin{lema} \label{regularidad hyperbolica}
The 2-form $\mu_N$ belongs to  $\esp\big(\XC{N}\big)$.
\end{lema}

\pruf It suffices to verify the assertion around the singularities of $\mu_N$, namely cusps and elliptic points. As $dd^c \log y = -\frac{dx \wedge dy}{4\pi y^2} = -\frac{\mu_1}{12}$, we only need to show that $\log y$ defines a $\Luno$ function on a neighborhood of these points. 

\underline{Cusps}: Let $s$ be the width of the given cusp, and let $\sg = \left(\begin{array}{cc}
* & * \\
u & v 
\end{array} \right) \in GL_2(\erre)$ be the element taking it to $\infty$. In terms of the local parameter $q=e^{2\pi i \sg (z)/s}$, we have the relation

\begin{equation}\label{puntas}
y= -(\log|q|)\frac{s}{2\pi}|j_\sg (z)|^2,
\end{equation}

\noindent where $j_\sg (z)= uz +v$. From this, it is easily checked that $\log y$ is $\Luno$ in a neighborhood of $q=0$. 

\underline{Elliptic points}: Let $z_0=x_0 + i y_0$ be an elliptic point of order $n$. A  local parameter is $\tau = \big( \frac{z-z_0}{z-\bar{z_0}}\big)^n$. After the choice of a branch of logarithm,  we have that 

\begin{equation} \label{elipticos}
 y = y_0 \frac{1-|\tau|^{2/n}}{|1-\tau^{1/n}|^2} \quad \textrm{almost everywhere}.
\end{equation}

 From this formula it is easy to verify that $\log y$ is $\Luno$ in a neighborhood of $\tau=0$. $\blacksquare$ \\

We introduce the notation $\Ch := \CH \big(\XZ{N}\big)_{\mu_N}$. 

\begin{thm}\label{Hecke aritmetico}

\begin{enumerate}

\item \label{Hecke bien definido}The correspondences $T_l$ with $l \nmid N$ and $w_d$ with $d | N, d \nmid N/d$,  induce homomorphisms $\hat{T}_l$ and $\hat{w}_d$ on $\CH (\XZ{N})$.

\item  \label{Hecke funtorial} These homomorphisms preserve $\Ch$.

\item \label{autoadjunto} The morphism $\hat{T}_l$ is self-adjoint. 

\item \label{Hecke conmutativo} Suppose $N$ squarefree. The sub-algebra of End$\Big(\CH \big(\XZ{N}\big)\Big)$ spanned by $\hat{T}_l$ with $l \nmid N$ and $\hat{w}_d$ with $d | N$  is commutative.

\end{enumerate}

\end{thm}

\pruf The assertion \eqref{Hecke bien definido} is an immediate  consequence of  Theorem \ref{main}, \eqref{bien definido}. 
It is easily verified that $T_l \mu_N = (l+1)\mu_N$. Then,  using  Lemma \ref{regularidad hyperbolica} and  Theorem \ref{main} \eqref{funtorial}, we conclude \eqref{Hecke funtorial}. Assertion \eqref{autoadjunto} follows from Corollary \ref{main en correspondencias}, \eqref{simetrico} and Lemma \ref{Heckesimetrico}.  
To prove \eqref{Hecke conmutativo}, we first note that the sub-algebra of End $\Big( Div \big(\XZ{N}\big) \Big)$ spanned by the correspondences  $T_l$ and  $w_d$ is commutative (cf. \cite{Shimura}, Proposition 3.32 to check commutativity for horizontal divisors and Lemma \ref{fibras finitas} for vertical divisors). Then we use Corollary \ref{commutacion aritmetica} to deduce the claim. $\blacksquare$ \\

 The next Lemma will be used to construct a $\mu_N$- admissible Green function for $D_\infty^N$.

\begin{lema} \label{formamodular}
\begin{enumerate}

\item \label{efe} Let $P_\infty:= \G(N) \infty \in \XC{N}$. There exists a weight $k \in 12 \ene$ cuspidal form $f \in S_k\big(\G(N)\big)$ and an integer $r \in \z_{>0}$ such that 
\begin{equation}\label{divisor}
\dvemph f = r [P_\infty] \textrm{ on }\XC{N}.
\end{equation}

 Moreover, $k/r = 12/ [\G(1): \G(N)]$.

\item \label{unicida} Suppose $f_1 \in S_{k_1} \big(\G(N)\big)$, with $k_1 \equiv 0 \mod 2$, satisfies the equality $\dvemph f_1 = r_1 [P_\infty]$ for some positive integer $r_1$. Then there exists $A \in \ce^*$ such that $$f^{r_1} = A f_1^r.$$ Moreover,  $k_1/r_1 = 12/ [\G(1): \G(N)]$.

\end{enumerate}
\end{lema}

\pruf  if $N=1$, we consider the  discriminant function

\[
\Delta(z) = q\prod_{n=1}^\infty (1-q^n)^{24} \in S_{12}\big(\G(1)\big), \quad q=e^{2\pi iz},
\]
\noindent  satisfying $\dv \Delta = [P_\infty]$. 

If  $N > 1$, we consider $\Delta$ as a modular form in $S_{12}\big(\G(N)\big)$.  We have that 

$$  \dv (\Delta)  =\sum_{C \textrm{ cusp}} a_C [C], $$

\noindent where $a_C > 0$ for every cusp $C$ and  $\sum_{C \textrm{ cusp}} a_C = [\G(1) : \G(N)] =:d_N$.  The divisor

\[
\dv (\Delta) - d_N[P_\infty]
\]

\noindent has degree 0 and is supported on the cusps. By the Manin-Drinfeld theorem (\cite{Drinfeld}, \cite{El}), there exists a meromorphic function $g$ on $\XC{N}$ and a positive integer $n$ such that

\[
\dv (g) = n (\dv (\Delta) - d_N[P_\infty]).
\]

Thus, $ f := \Delta^n / g$ is an holomorphic modular form in $S_{12n}\big( \G(N)\big)$  that satisfies \eqref{divisor}. We have that

\[
\dv (f) = n \dv \Delta - \dv g = n d_N [P_\infty],
\]

\noindent so $r = n [\G(1) : \G(N)]$ and $k/r = 12/ [\G(1) : \G(N)]$, which proves the first assertion. To prove the second, we remark that since the  holomorphic modular form $f_1^r/f^{r_1}$ has no zeros, it must be constant. We obtain that 
$k_1/r_1 = k/r = 12/ [\G(1) : \G(N)]$, i.e. by writing down that its weight is 0. $\blacksquare$ \\

Let $f$ be as in part \eqref{efe} of the preceding Lemma. We define

\begin{equation} \label{greensina}
g_\infty(z):= -\frac{1}{r} \log |f(z)^2 y^k|.
\end{equation}

The statement \eqref{unicida} in the Lemma implies that $g_\infty$ depends on the choice of $f$ only up to an additive  constant. An alternative proof of this fact is given by part \eqref{parteteorica} of the following proposition:

\begin{prop} \label{greenalinfinito}

\begin{enumerate}

\item \label{parteteorica} $g_\infty$ induces a $\mu_N$-admissible  $L^2_1$-Green function for  $[P_\infty]$.

\item  \label{expansion} For $c \in \p (\q)$, put $Q = \G(N) c$. We have that

\[
g_\infty (z) = -\delta_{P_\infty,Q}\log|q|^2 -b_N \log(-\log(|q|)) + b, \quad
z \rightarrow c,
\]

\noindent where $q$ is the standard local chart around $Q$, the function $b$ is  $C^\infty$   and $b_N = 12/[\G(1)
: \G(N)]$.
\end{enumerate}
\end{prop}

\pruf  let $z_0 = x_0 + i y_0 \in \h$ be a point of order
$n$ and let  $\tau$ be the local chart around $z_0$ as in 
\eqref{elipticos}.  If $n=1$ (i.e. if $z_0$ is not an elliptic point), the expression \eqref{elipticos}, shows  that 
 $g_\infty|_{\YC{N}}$ is $C^\infty$ at $z_0$. If $n>1$,  $g_\infty$ is
continuous at $z_0$ and belongs to $L^2_1$. We have that

\begin{eqnarray}
dd^c g_\infty (z) & = & -\frac{k}{r} dd^c \log y \nonumber \\
& = & - \frac{12}{[\G(1) : \G(N)]} \Big(-\frac{\mu_1}{12}\Big) \nonumber \\
& = & \mu_N. \nonumber
\end{eqnarray}

The $q$-expansion of $f$ at $Q$ is given by

\[
f|_k \sigma (z) = j_\sigma(z)^k f(z) = \sum_{n=r_Q}^\infty a_n q^n, \quad a_{r_Q}
\neq 0,
\]

\noindent where $r_Q $ is the order of vanishing of $f$ at $Q$.
Using the relation \eqref{puntas}, we obtain that

\[
y^k|f(z)|^2 = (-\log|q|)^k |q|^{2 r_Q} |h(q)|^2,
\]

\noindent where $h$ is an holomorphic function  defined on a disc containing 0 and such that $h(0) \neq 0$. We find that

\[
g_P(z) = -\delta_{P,Q} \log |q|^2 -\frac{k}{r_Q}\log(-\log|q|) -
\frac{1}{r_Q} \log |h(q)|^2
\]

\noindent and $b_N = k/r_Q$ by Lemma \ref{formamodular}, \eqref{efe}. $\blacksquare$ \\

Let $D_\infty:= D_\infty^N$.  We denote by $\hat{D}_\infty=(D_\infty, g_\infty)$ the  compactified divisor obtained by the proposition above.

\begin{lema} \label{Hecke infinito}
We have that $$\hat{T}_l \hat{D}_\infty = (l+1) \hat{D}_\infty + (0,c_{N,l}),$$  with $c_{N,l} = \frac{12(l-1)}{[\G (1) : \G (N)]} \log(l)$.
\end{lema}

\pruf as $T_l D_\infty = (l+1)D_\infty$ (cf. Lemma \ref{dinfinito}), the difference $T_l g_\infty - (l+1)g_\infty$ is a constant $c$. To evaluate $c$, we compare the expansion around $\infty$ of both functions.  With the notations of equation \eqref{greensina}, we consider the expansion at $\infty$ of $f$:

\[
f(z) = a_r q^r + a_{r+1} q^{r+1} + \ldots, \qquad a_r \neq 0, \quad r \geq 1.
\]

We  write $|f(z)|^2 = |a_r q^r|^2 h(q)$, where $h$ is defined on a neighborhood of 0 and $h(0)=1$. 
On the other hand, if we put 

\begin{equation} \label{matrices}
\alpha_j : = \left\{ \begin{array}{ll} \left(\begin{array}{cc}
1 & j \\
0 & l
\end{array}\right) & \textrm{ if } 0 \leq j \leq l-1 \\
 & \\
\left(\begin{array}{cc}
l & 0 \\
0 & 1
\end{array}\right) & \textrm{ if } j=l.
\end{array}\right.
\end{equation}

\noindent then we have that

\begin{equation} \label{eliminac2}
T_l g_\infty (z) = -\frac{1}{r}\sum_{j=0}^l \log |\big(\textrm{Im} (\alpha_j z)\big)^k f(\alpha_j z)^2|.
\end{equation}

But,

\[
\textrm{Im} (\alpha_j z) = \left\{ \begin{array}{ll}
-\frac{1}{2\pi l } \log(|q|) & \textrm{si } 0 \leq j \leq l-1 \\
-\frac{l}{2\pi} \log |q| & \textrm{si } j=l
\end{array}\right.
\]

\noindent and 

\[
|f(\alpha_j z)| = \left\{ \begin{array}{ll}
|a_r q^{r/l} | h_j (q) & \textrm{if } 0 \leq j \leq l-1 \\
|a_r q^{rl}| h_l(q) & \textrm{if } j=l,
\end{array}\right.
\]

\noindent where the functions $h_0, \ldots, h_l$ are defined around 0, positive, continuous and satisfy $h_j(0) = 1$.  Using these notations, we obtain that

$$T_l g_\infty(z)  - (l+1)g_\infty (z) = \frac{k}{r} (l-1)\log l - \frac{1}{r} \sum_{j=0}^{l} \log |h_j(q)|^2 + \frac{1}{r} \log |h(q)|^2.$$

Hence, taking the limit as $z \rightarrow \infty$ and using Lemma \ref{formamodular}, \eqref{unicida} we obtain that $c = b_N \log(l^{l-1}). \blacksquare$ \\

\begin{definiciones}  \label{i infinito}

\begin{itemize}

\item For a divisor $D \in $Div$(\XC{N}_{\q})$, we denote by $\overline{D}$ its Zariski closure in $\XZ{N}$. If $\deg D = 0$, we denote by $\Phi (D)$ a vertical divisor (which may have rational coefficients) such that $\overline{D} +\Phi (D)$ has degree zero on every irreducible component of every vertical fiber of $\XZ{N}$. The divisor $\Phi (D)$ is well defined up to  a sum of rational multiples of finite fibers.

\item Suppose that the genus of $\XC{N}$ is nonzero. We identify $\JNQ(\q)$ with $Pic^0 (\XQ{N})$ using the rational point $\G(N) \infty$. Let $$i_\infty : \JNQ(\q) \rightarrow \Ch$$ be given by $i_\infty(D) = (\overline{D} + \Phi(D), g_D )$, where $g_D$ is a $\mu_N$-admissible $\Luno$-Green function for $D$ normalized by the condition $$ \izq  i_\infty(D) , \hat{D}_\infty \der = 0.$$ As $\deg D=0$, the  function $g_D$ does not depend on the choice of a normalization for $g_\infty$ (cf. example \ref{grado cero}).

\end{itemize}
\end{definiciones}

\begin{lema}
The function $i_\infty$ is  well defined. Moreover, it is an embedding of abelian groups.
\end{lema}

\pruf We must check that $i_\infty(D)$ does not depend of the choice of $\Phi(D)$. Let $X_p:= \XZ{N} \otimes \Fp$ be a finite fiber. Let $g_p$ be a   $\mu_N$-admissible $\Luno$-Green function for $D$ such that $$ \izq (\overline{D} +\Phi (D) + X_p, g_p) , \hat{D}_\infty \der = 0.$$ 

As $(X_p, g_p - g_D)= (\overline{D} +\Phi (D) + X_p, , g_p)-i_\infty(D)$ is a compactified divisor, the difference $g_p - g_D$  is a constant $c$. Moreover, using the example \ref{grado cero},  the equality $$\izq (X_p, c), \hat{D}_\infty \der = 0$$ unwinds to $$\log p + \frac{c}{2} = 0.$$ This gives $(X_p, g_p - g_D) = (\dv p, -\log p^2)$, that is a (compactified) principal divisor. 

The injectivity of $i_\infty$  then follows plainly from the definitions. $\blacksquare$\\

\begin{cor} \label{compatibilidad}

\begin{enumerate}

\item \label{mismos Heckes} The following diagram commutes:

\[
\xymatrix{
 \JNQ(\q) \ar[r]^{T_l} \ar@{^{(}->}[d]^{i_\infty} & \JNQ(\q) \ar@{^{(}->}[d]^{i_\infty} \\
 \Ch \ar[r]^{\hat{T}_l} & \Ch. }
\]

\item Let $\izq , \der_{NT}$ denote the Néron-Tate height pairing on $\JNQ(\q)$. We have that

$$  \izq i_\infty x,  i_\infty y \der =- \izq x , y\der_{NT},$$ for all $x,y \in \JNQ(\q)$.
 
\end{enumerate}
 \end{cor}

\pruf The second assertion is just a restatement of the Faltings-Hriljac formula (\cite{Faltings}, \cite{Hriljac}, \cite{Moret-Bailly} 6.15).
Using the notations of the definitions \ref{i infinito}, the first assertion amounts to checking that, for every divisor $D$ on $\XZ{N}$, we have the equality
\begin{equation}\label{conmutacion}
T_l g_D = g_{T_l D}.
\end{equation}

Theorem \ref{Hecke aritmetico}, \eqref{Hecke funtorial} ensures that  $T_l g_D $ is an admissible $\Luno$-Green function for $T_l D$. Hence,   \eqref{conmutacion} holds up to a constant. The constant must be zero since by Theorem \ref{Hecke aritmetico} \eqref{autoadjunto} and Lemma \ref{Hecke infinito},  we have that 

\begin{align}
\izq (T_l D, T_l g_D) , \infinito \der &= \izq (D, g_D) , \hat{T}_l \infinito \der & \nonumber \\
 &= (l+1) \izq (D,g_D) , \infinito \der + \izq (D,g_D) , (0,c_{N,l}) \der & \nonumber \\
 &= 0, & \nonumber
\end{align}

\noindent i.e., by the definition of $g_D$ and the fact that $\deg D =0$ (cf. Example \ref{grado cero}).  $\blacksquare$ \\

Theorem \ref{Hecke aritmetico}, \eqref{autoadjunto} and the previous Corollary give the following result:

\begin{cor}
The Hecke correspondence $T_l : \JNQ(\q) \rightarrow \JNQ(\q)$, with $l \nmid N$, is
self-adjoint with respect to the N{\'e}ron-Tate height pairing.
\end{cor}

\begin{rem}
This result is not new (cf. \cite{GrossZagier}). It can be proved  using the fact that $\JN(\q)\otimes \ce$ and $S_2(\G(N))$ are isomorphic as Hecke modules. The latter algebra is diagonalizable,  which implies that the former is also diagonalizable, which in turn implies the self-adjointness.  
\end{rem}

\subsection{Arithmetic Chow group of $\XZ{N}$ with real coefficients} \label{num}
In this Section, $N$ is a fixed squarefree integer. Let 
$\widehat{Div}_{\erre} (N)$ denote the $\erre$-vector space made of pairs $(D,g)$, with $D$ a Weil divisor with real coefficients on $\XZ{N}$ and $g$ a $\mu_N$-admissible  $L^2_1$-Green function for $D(\ce)$. Let  $P_\erre(N)$ be the sub-espace spanned by the principal divisors. We put

\[
\ChR := \widehat{Div}_{\erre} (N) / P_\erre(N).
\]

\begin{definicion} Let

\[
K_N  = \{ x \in \ChR | \mathopen<x,y\mathclose> = 0, \forall y \in \ChR\}.
\]
We define the arithmetic Chow group up to numerical equivalence as

\[
\Chnum := \ChR /K_N
\]
\end{definicion}

\begin{rems}
\begin{itemize}
\item The arithmetic intersection pairing $<, >$ extends to a nondegenerate bilinear form on $\Chnum \times \Chnum$.

\item Let $\sg : \Ch \rightarrow \Chnum$ be the natural map. We have an exact sequence

\[
\xymatrix{
0  \ar[r] & \JNQ(\q)_{tors}  \ar[r]^-{i_\infty}  & \Ch \ar[r]^-{\sg} & \Chnum. }
\]

\noindent  Indeed, if $\izq \sg(x), y\der =0 $ for every $y \in \ChR$, then by varying $y$ over all classes of compactified divisors of the form $(F,0)$ and $(0,c)$ with $F$ an irreducible component of a fiber of $\XZ{N}$ and $c\in \erre$ we conclude that the underlying divisor of any representative of $x$ comes from a divisor $\tilde{x} \in \JNQ(\q)$. Then the  Néron-Tate height of    $\tilde{x} $ vanishes because of the Faltings-Hriljac formula, implying that it is a torsion element.

\end{itemize}
\end{rems}

For $u \in \{ 0, \infty \}$ and $p|N$, we put  $\hat{X}_p^u := (X_p^u(N),0)$. 
Let $\hat{G}_p := \hat{X}_p^\infty - \hat{X}_p^0$ and  let $F := \{ (0,c), c \in \erre\}$. We define

$$\Eis := (F \oplus \erre \hat{D}_\infty) \oplus
\big(\overset{\perp}{\oplus}_{p|N} \erre \hat{G}_p\big) \subset \Chnum.$$

\noindent Here, the notation $A \overset{\perp}{\oplus} B$ means that $A \oplus B$ is a direct sum of vector subspaces and that $A$ and $B$ are mutually orthogonal. We remark that the space $F \oplus \erre \hat{D}_\infty$ does not depend on the normalization of $g_\infty$.

The Hodge index theorem in this context can be written as

\begin{thm}\label{modulos de Hecke}
We have that

\[
\Chnum = \Eis \overset{\perp}{\oplus} J,
\]

\noindent where $J$ is, by definition, the orthogonal complement of $\Eis$. Moreover, the rule $$(D,g) \mapsto D$$ induces an isomorphism of Hecke modules $J \cong \JNQ (\q) \otimes \erre$. 

\end{thm}

The proof of the Hodge index theorem given in \cite{Moret-Bailly}, p. 85, works in this situation (cf. also \cite{Bost}, Theorem 5.5). We just remark that the choice of basis $\hat{G}_p$ is convenient because it forces an element in the orthogonal complement of $\Eis$ to have degree 0 on every irreducible component of every finite fiber. The compatibility of the Hecke actions on $J$ and $ \JNQ (\q) \otimes \erre$ is given by Corollary \ref{compatibilidad} \eqref{mismos Heckes}.

The following assertion is a simple consequence of Lemma \ref{fibras finitas} and Lemma \ref{Hecke infinito}.

\begin{prop}
 
\begin{enumerate}

\item \label{facilona} The spaces $F$ and $\erre \hat{G}_p$
 with $p | N$ are eigenspaces for $\hat{T}_l$ and   $\hat{w}_d$ with $d | N$. More precisely, $\hat{T}_l x = (l+1)x $ for all $x \in F \oplus   \hat{G}_p$ and the morphism  $\hat{w}_d$ is the identity on $F$.  On $\hat{G}_p$, $\hat{w}_d$ is the identity (resp. $-\hat{w}_d$ is the identity)  if $w_d(X_p^\infty)= X_p^\infty$ (resp. if $w_d(X_p^\infty)= X_p^0$). In particular, $w_N|_{\hat{G}_p}=-id_{\hat{G}_p}$.

\item \label{nodiagonalizable} The space
$F \oplus \erre \hat{D}_\infty$ is stable under $\hat{T}_l$. More precisely,

\[
\hat{T}_l (0, c) = (l+1) (0,c), \quad \hat{T}_l \hat{D}_\infty = (l+1)\hat{D}_\infty + (0,c_{N,l}),
\]

\noindent with $c_{N,l} = \frac{12(l-1)}{[\G (1) : \G (N)]} \log(l)$.

\end{enumerate}

\end{prop}

\subsection{Self-intersection and refined self-intersection of the dualising sheaf} \label{explicacion de la normalizacion}

Let $\omega$ be the dualising sheaf on $\XZ{N}$. This sheaf induces a class in $CH^1(\XZ{N})$, and  we also denote by $\omega$ a divisor in this class. Then
$\omega_{\q} := \omega \otimes \q$ is a canonical divisor on $\XQ{N}$.  Let $\pi : \XQ{N} \rightarrow \XQ{1}$ be the natural morphism and let $\bf{i}$ (resp. $\bf{j}$) be the orbit of $\sqrt{-1}$ (resp. $e^{i\pi/3}$) in $\XC{1}$.  Using Hurwitz's formula and the Manin-Drinfeld  theorem (\cite{El}, \cite{Drinfeld}) we find that

\begin{equation} \label{CanonicoHeegner}
\omega_{\q} = (2g-2) [\infty]   -  H_i -
2H_j \qquad \textrm{ in } CH^1 (\XQ{N}) \otimes \q.
\end{equation}

Here, $g$ is the genus of $\XC{N}$, and the divisor $H_i$ (resp. $H_j$) is the sum of all divisors of the form $\frac{1}{2}\big([P]-[\infty]\big)$ (resp. $\frac{1}{3} \big( [P]-[\infty]\big)$)  with $\pi(P) = \bf{i}$ (resp.   $\pi(P) = \bf{j}$) and $\pi$ unramified at $P$. Equivalently, the points $P$ appearing in $H_i$ (resp. $H_j$) are the Heegner points of discriminant -4 (resp. -3) on $\XC{N}$ (cf. \cite{MichelUllmo}, Section 6). 

We introduce the following notation. For a divisor $E \in $Div$(\XQ{N})$, we write $D_E$ for its Zariski closure in $\XZ{N}$. We deduce from (\ref{CanonicoHeegner}) the equality

\begin{equation} \label{DualizanteHeegner}
\omega =  (2g-2) D_\infty  - D_{H_i} -
2 D_{H_j} + V \qquad \textrm{ in } CH^1\big(\XZ{N}\big)\otimes\q,
\end{equation}

\noindent where $V$ is a vertical divisor contained inside the space of fibers of bad reduction. This comes from the fact that both sides of (\ref{DualizanteHeegner}) are equal on the generic fiber.

Let $W$ be a vertical divisor such that  $\omega_J:= \omega - (2g -2) D_\infty + W$ has degree zero on each irreducible component of every fiber of $\XZ{N}$. The divisor $\omega_J$ is then  identified with a  point in $\JN(\q)$.  We define $\wbar_J := i_\infty(\omega_J) \in J$ (where $\wbar_J$ does not depend on the choice of $W$) and

\[
\omega_{\Eis} := \omega - \omega_J. 
\]

We will now specify a compactification of $\omega_{\Eis}$. We define $\hat{W}=(W,c)$, with $c$ a constant such that $\izq \hat{W}, \hat{D}_\infty \der = 0$. To normalize the function $g_\infty$ underlying $ \hat{D}_\infty$, we present the following considerations: let $C$ be a divisor supported on the cusps with the additional property  that $D_\infty$ is not in its support. Set $D= D_\infty + C$, and let $g$ be a $\Luno$-admissible Green function for $D$. Proposition \ref{greenalinfinito} ensures that $$a(g):=  \lim_{z \rightarrow i \infty} g (z) + \log|q|^2 + b_N \log(-\log |q|)$$ is a well defined real number. We normalize $g_\infty$ by imposing $a(g_\infty) = -\frac{12\log (2)}{[\G(1) : \G(N)]^2}$.  

This particular choice is motivated by the theory of metrized line bundles. Let $\overline{\mathcal{M}} = (\mathcal{M}_{12} (\G(N)), \Vert \cdot \Vert_{Pet})$ be the line bundle of weight 12 modular forms on $\XZ{N}$ endowed with the  Petersson metric, as defined in 
\cite{Kuhn}. Namely, if $f$ is a section of $\mathcal{M}_{12} (\G(N))$, then

\begin{equation} \label{nn} 
\Vert f_{\ce}(z)\Vert_{Pet}^2 = |f_{\ce}(z)|^2 (4 \pi y )^{12}, \quad z = x+iy \in \h.
\end{equation}

This  factor $4 \pi$ in \eqref{nn} is motivated by a natural isomorphism between the line bundle of weight 12 modular forms attached to $\Gamma(N)$  and the 12th power of the line bundle on $X(N)$ induced by the canonical bundle on the  universal elliptic curve. This last line bundle has a canonical metric (the so-called $L^2$ metric) and $\Vert \cdot \Vert_{Pet}$ in this case is obtained by following the isomorphism (see \cite{Kuhn}, Section 4.14).

We denote by $\Delta$ the section of $\mathcal{M}_{12} (\G(N))$ defined by the discriminant function. Using Proposition  \ref{greenalinfinito}, we see that the pair
\[
\hat{E} := (\dv (\Delta), -\log \Vert\Delta \Vert_{Pet}^2)
\]

\noindent is a compactified divisor. The product expansion of the discriminant function $\Delta$ shows that $$ a(-\log \Vert\Delta \Vert_{Pet}^2) =  -\frac{12\log(2)}{[\G(1) : \G(N)]}.$$

 This data defines an element  $\wbar_{\Eis} \in \Eis$. Finally, we have the element  $\wbar = \wbar_{\Eis} + \wbar_J \in \Chnum$. The main ingredient needed for  the computation of $\wbar_{\Eis}^2$ is  a calculation due to U. Kühn. We  deal mostly with the intersections at the bad fibers. Let us postpone the proof of Theorem \ref{calculo Eis} to the following Section.

If there are no elliptic points on $\XC{N}$ (i.e. if there are two different prime divisors $p,q |N$ such that $p, q \notin \{2,3\}$,  $p \equiv 3$ mod 4 and $q \equiv 2 $ mod 3), then there are no nonramified points above
the elliptic points of $\XC{1}$ (\cite{Shimura}, Proposition 1.43). In this case, $\wbar_J = 0$ because the divisors 
$D_{H_i^0}$ et $D_{H_\rho^0}$  in \eqref{DualizanteHeegner} are empty. \\

Theorem \ref{modulos de Hecke} allows to associate to a normalized eigenform $f \in S_2\big(\G(N)\big
)$ the $f$-isotypical component of $\wbar_J$, which we denote by $\wbar_f$. The Gross-Zagier formula and the Gross-Kohnen-Zagier theorem then allow us to compute the self-intersection of this element in some cases.\\

\noindent {\bf Proof} of Theorem \ref{calculo isotipico}: since both $\q (i)$ and $\q (\sqrt{-3})$ have class number one, the Heegner points in $H_i$ and $H_j$ are defined over these respective quadratic fields. Moreover, since Gal$(\q(i)/\q)$ preserves the set of Heegner points of discriminant $-4$, the divisor $H_i$  is defined over $\q$. The same argument applies to show that $H_j$ is also defined over $\q$.

Let $X^*$ be the quotient of $\XC{N}$ by the Fricke involution $w_N$, and let $J^*$ be the jacobian of $X^*$. Let $ \kappa : \XC{N} \rightarrow X^*$ be the canonical map, which is of degree two.  Let $\overline{\omega}_f := \kappa_* (\omega_f)$.

We remark  that we have the equalities 

\begin{eqnarray} 
h_{NT} \big( \overline{(H_i)}_f \big) = 2 h_{NT}  \big( (H_i)_f \big), & \quad h_{NT} \big( \overline{(H_j)}_f \big) = 2 h_{NT}  \big( (H_j)_f \big), \label{pasajes}\\
h_{NT} \big(\overline{\omega}_f \big) = 2 h_{NT}  \big( \overline{\omega}_f \big). & \nonumber
\end{eqnarray}

Indeed, we   have that

\begin{eqnarray}
2 h_{NT}\big(\overline{\omega}_f\big)  &=&  \izq \overline{\omega_f}, \kappa_* \kappa^* \overline{\omega_f} \der \nonumber \\
 &=&\izq  \kappa^* \overline{\omega_f}, \kappa^*  \overline{\omega_f} \der \nonumber \\
 &=& h_{NT} \big(  \omega_f +  (w_N)_* \omega_f \big) \nonumber \\
 &=&  h_{NT} \big(  2 \omega_f \big) \nonumber \\
 &=& 4 h_{NT} \big(  \omega_f \big), \nonumber
\end{eqnarray}

\noindent where the last equality comes from the fact that the Heegner points of given discriminant are permuted under $w_N$. The same argument applies to check the other two equalities in \eqref{pasajes}. We then have that

\begin{eqnarray}
-2\wbar_f^2 &=& 2h_{NT} \big(\omega_f\big)\nonumber \\
&=& h_{NT}\big(\overline{\omega}_f\big) \nonumber \\
 &=&  h_{NT} \big(\overline{(H_i)}_f\big) + 4 h_{NT}  \big(\overline{(H_j)}_f\big) + 4 \izq \overline{(H_i)}_f , \overline{(H_j)}_f \der_{NT} \nonumber
\end{eqnarray}

It follows from the Gross-Kohnen-Zagier theorem \cite{GKZ} that the images in 

\noindent $\big(\JN(\q)/w_N\big)_f\otimes \erre$ of the divisors $H_i, H_j$ are collinear (see ch. V of loc. cit. for the case of even discriminant). In particular,  $$\izq \overline{(H_i)}_f , \overline{(H_j)}_f \der_{NT} = h_{NT}\big(\overline{(H_i)}_f\big)^{1/2}  h_{NT}\big(\overline{(H_j)}_f\big)^{1/2}.$$ This, together with the equalities \eqref{pasajes}, prove \eqref{alturas}. 

Suppose that $N$ is prime. Then $H_i$ is either the empty divisor, or a divisor that contains two Heegner points that are interchanged by $w_N$. In this situation, the calculation of $h_{NT} (H_i)_f$ in terms of special values  is given by \cite{BY}, Corollary 7.8. The same argument applies to $H_j$ and the calculation of $h_{NT} (H_j)_f. \blacksquare$ 

\begin{rems}

\begin{itemize}

\item If $N=p_1 p_2 \ldots p_s$ is a squarefree integer with $s>1$, then $(H_j)_f$ is either empty or contains $2^s$ ($f$-isotypical components of) Heegner points. In the latter  case, they give raise to $2^{s-1}$ (not necessarily different) elements of  $\Big(\JN\big(\q(\sqrt{-3})\big)_f/w_N \Big)\otimes \erre$. The Gross-Kohnen-Zagier theorem and the Gross-Zagier formula ensure that they are all collinear, and that they all have the same height. However, we do not know how to exclude the possibility that any given two of these points have opposite signs.  This prevents us from computing $h_{NT} \big( ( H_j)_f \big)$. In the case that $s=1$, the technical improvements on the Gross-Zagier formula by Bruinier and Yang rule out this problem. 

\item We also point out that the original Gross-Zagier formula as stated in \cite{GrossZagier} does not apply to Heegner points of even discriminant, yet their methods undoubtedly apply in the case needed here, namely discriminant -4. On the other hand, the restriction on the parity has been eliminated in recent works (cf. \cite{BY}, \cite{Z}, \cite{CM}).  

\end{itemize}

\end{rems}

\subsubsection{Computation of $\wbar_{\Eis}^2$} This Section is devoted to the proof of Theorem \ref{calculo Eis}. As  $\wbar_{\Eis}^2=  (2g - 2)^2 \hat{D}_\infty^2 + \hat{W}^2$, it suffices to prove the following two lemmas:

\begin{lema}  (U. K\"uhn, \cite{Kuhn}) We have that $$\hat{D}_\infty^2  =  \frac{144}{[\G(1) : \G(N)]}\Big(\frac{1}{2}\zeta(-1) + \zeta'(-1)\Big). $$

\end{lema}

\pruf It is shown in \cite{Kuhn}, Corollary 6.2, that the generalized self-intersection number $\overline{\mathcal{M}}^2$ is given by

\[
\overline{\mathcal{M}}^2 = 144 [\G(1) : \G(N)]\Big(\frac{1}{2}\zeta(-1) + \zeta'(-1)\Big).
\]

 By   \cite{Kuhn} Proposition 7.4,  we have that  $\hat{E}^2 = \overline{\mathcal{M}}^2$. The  Manin-Drinfeld Theorem (\cite{Drinfeld}, \cite{El}) implies that

\[
\dv \Delta = [\G(1) : \G(N)] D_\infty + V \quad \textrm{ in } CH^1(\XZ{N}) \otimes \q,
\]

\noindent where $V$ is a vertical divisor. Thus, $$ [\G(1) : \G(N)] \hat{D}_\infty = \hat{E} + (0,c) \quad \textrm{ in } \Chnum.$$

As $a(g_\infty) =a( -\log \Vert \Delta\Vert_{Pet}^2)$, by definition, we obtain that $c=0$. We conclude that

\[
\hat{D}_\infty^2 = [\G(1) : \G(N)]^{-2}\hat{E}^2,
\]

\noindent as required to conclude the proof. $\blacksquare$ \\

\begin{lema} \label{wcuadrado}
We have that $$\hat{W}^2  =  -(g-1)^2 \sum_{p|N}  \frac{\log p}{g - 2g_p +1}.$$
\end{lema}

We break the proof of this lemma in two further lemmas.

\begin{lema} \label{intersecciones}
For $p \nmid N$, we have that

\[
\mathopen<W , X_p^\infty\mathclose> = (g-1)\log p, \quad \mathopen<W , X_p^0\mathclose> = (1-g)\log p.
\]

Moreover, we have that

\begin{equation} \label{componentes}
\mathopen<X_p^\infty , X_p^0\mathclose> = (g - 2g_p + 1) \log p.
\end{equation}
\end{lema}

\pruf : We have  $\dv p = (X_p^\infty + X_p^0, -\log p^2)$, for all  $p | N$. Thus,

\begin{eqnarray}
0 & = & \mathopen< \hat{X}_p^\infty , (X_p^\infty + X_p^0, -\log p^2) \mathclose> \nonumber \\
& = & \mathopen< X_p^\infty , X_p^\infty + X_p^0\mathclose> \nonumber \\
& = & (X_p^\infty)^2 + \mathopen<X_p^\infty ,  X_p^0\mathclose>. \nonumber
\end{eqnarray}

Similarly,  $ \mathopen<X_p^0 , X_p^\infty\mathclose> + (X_p^0)^2 = 0.$ Let 

\[
W=\sum_{p|N} W_p, \quad |W_p| \subseteq \XZ{N} \otimes \Fp.
\]

By the adjunction formula (\cite{Liu}, Chapter 9, Theorem 1.37):

\begin{eqnarray}
(2 g_p - 2) \log p & = & \mathopen<\omega , X_p^\infty\mathclose> + (X_p^\infty)^2 \nonumber \\
& = & \mathopen<\omega_{\Eis} , X_p^\infty \mathclose> + (X_p^\infty)^2 \nonumber \\
& = & (2g-2) \mathopen<D_\infty , X_p^\infty\mathclose> - \mathopen<W_p , X_p^\infty\mathclose> + (X_p^\infty)^2 \nonumber \\
& = & (2g-2) \mathopen<D_\infty , X_p^\infty\mathclose> - \mathopen<W_p , X_p^\infty\mathclose> - \mathopen<X_p^\infty , X_p^0\mathclose>.
\end{eqnarray}

For the same reasons, we have that
\[
(2 g_p - 2) \log p = (2g-2) \mathopen<D_\infty , X_p^0\mathclose> - \mathopen<W_p , X_p^0\mathclose> - \mathopen<X_p^\infty , X_p^0\mathclose>.
\]

As $\mathopen<D_\infty , X_p^\infty\mathclose> = \log p$ and $\mathopen<D_\infty , X_p^0\mathclose> = 0$, we obtain that

\begin{eqnarray}
\mathopen<W_p , X_p^\infty\mathclose> & = & 2(\log p)(g-g_p) - \mathopen<X_p^\infty , X_p^0\mathclose>, \label{remate} \\
\mathopen<W_p , X_p^0\mathclose> & = & -(\log p )(2g_p - 2)- \mathopen<X_p^\infty , X_p^0\mathclose>.\nonumber
\end{eqnarray}

To cumpute $\mathopen<X_p^\infty , X_p^0\mathclose>$, we use the equality

\begin{eqnarray}
0 & = & \mathopen<(W_p,0) ,  (X_p^\infty + X_p^0, -\log p^2) \mathclose> \nonumber \\
& = & \mathopen<W_p , X_p^\infty + X_p^0\mathclose> \nonumber \\
& = & \mathopen<W_p , X_p^\infty\mathclose>+ \mathopen<W_p , X_p^0\mathclose>. \nonumber
\end{eqnarray}

Using (\ref{remate}), we obtain equation \eqref{componentes}. Using this result in  \eqref{remate} then finishes the proof. $\blacksquare$ \\

\begin{lema}\label{constanteinutil}

For some $c \in \erre$, we have the equality 

\[
\hat{W} = -\frac{(g-1)}{2}\sum_{p|N} (g-2 g_p +1)^{-1} \hat{G}_p + (0,c) \quad \textrm{ \emph{in} }\Chnum.
\]

\end{lema}

\pruf:  We write $\hat{W} = \sum_{p|N} a_p \hat{G}_p + (0,c)$ in $\Chnum$. Using Lemma \ref{intersecciones}, we find that

\begin{eqnarray}
(g-1)\log p & = & \mathopen<W , X_p^\infty\mathclose>  \nonumber \\
& = & \mathopen<\hat{W} , \hat{X}_p^\infty\mathclose>  \nonumber \\
& = & a_p \mathopen<G_p , X_p^\infty\mathclose> \nonumber \\
& = & a_p \big((X_p^\infty)^2 - \mathopen<X_p^0 , X_p^\infty\mathclose> \big) \nonumber \\
& = & -2a_p \mathopen<X_p^0 , X_p^\infty\mathclose> \nonumber \\
& = & -2a_p (g- 2g_p +1) \log p, \nonumber
\end{eqnarray}

\noindent where we have used \eqref{componentes} in the last equality. $\blacksquare$ \\

\noindent {\bf Proof } of Lemma \ref{wcuadrado}: using Lemma \ref{constanteinutil}, we have that

\[
\hat{W}^2 = \Big(\frac{g-1}{2}\Big)^2 \sum_{p|N} (g-2 g_p +1)^{-2} \hat{G}_p^2
\]

Using equation \eqref{componentes}, we obtain that
\begin{eqnarray}
\hat{G}_p^2 & = & \big((X_p^\infty)^2 + (X_p^0)^2 - 2 \mathopen<X_p^\infty , X_p^0\mathclose> \big)\nonumber \\
& = & -4 \mathopen<X_p^\infty , X_p^0\mathclose>  \nonumber \\
& = & -4(g- 2g_p +1) \log p, \nonumber
\end{eqnarray}

\noindent which concludes the proof. $\blacksquare$ \\

\section{Appendix: cohomology of $\Luno$ spaces on compact Riemann surfaces}

The goal of this Section is to prove Lemma \ref{Poisson} used in Section \ref{Chowgroup} to prove the existence of admissible $\Luno$-Green functions (Proposition \ref{existencia}). The strategy to solve equation \eqref{ecuaciondePoisson} is to first solve it locally, then to patch together the solutions to obtain a global one using cohomological machinery. The key local input to make this machinery work is the following 

\begin{thm} \label{local}
Let $U \subset \ce$ be a connected open set. Then  for every $g, f_1, f_2 \in L^2(U)$ the equation $$\Big( \frac{\partial^2}{\partial x^2} +  \frac{\partial^2}{\partial y^2} \Big) u = g + \frac{\partial f_1}{\partial x} + \frac{\partial f_2}{\partial y}$$ has a solution $u \in \Luno(U)$.

\end{thm}

This result is well known and can be handled by standard methods of PDE's (e.g. the Lax-Milgram Theorem and the Fredholm alternative). A proof of this theorem can be found in \cite{GilbargTrudinger}, p. 170, Theorem 8.3. We will also need an appropriate version of the $\bar{\partial}$ lemma:

\begin{lema} \label{deltalema}
Let $U \subset \ce$ be a connected open set, and let $f \in L^2(U)$. There exists a $u \in \Luno(U)$ such that $$\bar{\partial} u = f d\bar{z}.$$

\end{lema}

\pruf  Theorem \ref{local} implies that there exists a $u_0 \in \Luno(U)$ such that $\partial \overline{\partial} u_0 = \partial f \wedge d\bar{z}$. Let $h:= \overline{\partial} u_0 - f d\bar{z}$. As  $\partial h =0$, we have that  $h$ is harmonic. By the ellipticity of $\bar{\partial}\partial $, we deduce  that $h$ is $C^\infty$ (and in fact anti holomorphic because it is annihilated by $\partial$). Then, by a special case of the Dolbeaut lemma (cf. \cite{Forster}, Theorem 13.2) there exists a $u_1 \in C^\infty(U)$ such that $ \overline{\partial} u_1 = h$. Then $u := u_0-u_1$ belongs to $\Luno(U)$ and verifies $\overline{\partial} u = f d\bar{z}. \blacksquare$

Let $X$ be a compact Riemann surface. 
We will use the following property of the cohomology on $X$:

\begin{thm} \label{coho}
Consider the following exact sequence of sheafs on $X$: $$0 \rightarrow  A \rightarrow  B \rightarrow^{\alpha}  C \rightarrow  0.$$ If $H^1(X,B)=0$, then $H^1(X,A) \cong  C(X)/ \alpha\big(B(X)\big)$.
\end{thm}

A proof of this theorem can be found in  \cite{Forster}, Theorem 15.13.\\

Let $\espa^{(0,1)}(X)$ (resp. $\espa^{(1,0)}(X)$) be the space of 1-currents on $X$ of type $(0,1)$ (resp. $(1,0)$) which are locally $L^2$. More precisely, an element $\omega \in  \espa^{(0,1)}(X)$ can be written on sufficiently small open sets  $U$ as $f d\bar{z}$ with $f\in L^2(U)$ (and similarly for elements in $\espa^{(1,0)}(X)$).  Out of this definition we construct sheafs $\espa^{(0,1)}$ and $\espa^{(1,0)}$ in the obvious way. If we consider $\Luno$ as a sheaf on $X$ we have  well defined morphisms $\bar{\partial} : \Luno \rightarrow \espa^{(0,1)}$ and $\partial  : \Luno \rightarrow \espa^{(1,0)}$.

\begin{lema}\label{H1}
We have that \begin{itemize}

\item $H^1(X, \Luno)=0$

\item $H^1(X,\espa^{(1,0)}) =0$.

\end{itemize}
\end{lema}

As the sheafs $\Luno$ and $\espa^{(1,0)}$ are stable under multiplication by $C^\infty$ functions, this lemma can be proved carrying over the argument of \cite{Forster}, Theorem 12.6, word by word with either of the sheafs $\Luno$ or  $\espa^{(1,0)}$. 

\begin{lema} \label{suma}
We have that $$\espa^{(0,1)}(X)=\bar{\partial}\Luno(X)\oplus \overline{\Omega}(X),$$
where  $\overline{\Omega}$ is the sheaf of anti holomorphic 1-forms on $X$. 
\end{lema}

\pruf The space $\espa^{(0,1)}(X)$ is endowed with the inner product $$(\omega, \eta) = -i \int_{X}^{} \omega \wedge \overline{\eta}.$$  For $f \in \Luno (X)$ and $\omega \in \overline{\Omega}(X)$, we have that $$(\omega, \overline{\partial} f ) = -i \int_{X}^{} d( \omega \overline{f} ) = 0.$$ This shows that the subspaces $\bar{\partial}\Luno(X), \overline{\Omega}(X)$ are orthogonal. Hence, it suffices to show that the codimension of the former equals the dimension of $ \overline{\Omega}(X)$, that is the genus of $X$.

Let $O$ be the sheaf of holomorphic functions on $X$. Consider the sequence

$$0 \rightarrow O  \rightarrow \Luno \rightarrow^{\bar{\partial}} \espa^{(0,1)} \rightarrow 0.$$

We claim that this sequence is exact. To check exactness in the middle, let $U\subset X$ be an open set, and let $f \in \Luno(U)$ be such that $\bar{\partial}f=0$. This implies that $f$ is harmonic, and so  it must be  $C^\infty$. But then, $f$ must be holomorphic because it is annihilated by $\overline{\partial}$. Since the surjectivity of $\Luno \rightarrow^{\bar{\partial}} \espa^{(0,1)}$ needs to be checked only locally, this map is surjective by Lemma \ref{deltalema}.

Using Lemma \ref{H1} and Theorem \ref{coho} we conclude that  $$H^1(X,O) \cong \espa^{(0,1)} (X) / \overline{\partial}\Luno  (X).$$ As the dimension of $H^1(X,O)$ equals the genus of $X$, this finishes the proof. $ \blacksquare$ \\

\noindent \textbf{Proof of Lemma \ref{Poisson}}: the  vanishing of the integral of $\mu$ is necessary, as can be checked using Lemma \ref{identidadutilisima}. To prove the reciprocal, we  claim

\begin{enumerate}

\item \label{palta} $d  \espa^{1,0}(X)=d  \espa^{0,1}(X) = \partial \overline{\partial} \Luno(X)$

\item \label{exacta} we have an exact sequence $$0 \rightarrow  \Omega \rightarrow \espa^{(1,0)} \rightarrow^{d} \esp \rightarrow 0. $$

\end{enumerate}

Claim \eqref{palta} follows from Lemma \ref{suma}. To check exactness in the middle in claim \eqref{exacta}, we just remark that $\overline{\partial} \partial  f =0 $ and that $f \in \Luno(U)$ implies $f \in C^\infty(U)$ because of the ellipticity of $\overline{\partial} \partial$. Then, $\partial f$ is holomorphic because it is annihilated by $\overline{\partial}$. Exactness on the right follows from Theorem \ref{local}.

 Using the claims \eqref{palta} and \eqref{exacta}, Lemma \ref{H1} and Theorem \ref{coho}, we have that $$H^1(X, \Omega) \cong \esp(X)/d\big(\espa^{(1,0)}(X)\big) = \esp(X)/\partial \overline{\partial} \Luno(X).$$ As the space $H^1(X,\Omega)$ is 1-dimensional (a consequence of Serre duality, cf. \cite{Forster}, Theorem 17.11), the image of $\mu$ in $\esp(X)/\partial \overline{\partial} \Luno(X)$ must vanish (since otherwise all elements in $\esp(X)$ would have integral equal to zero). This proves  the existence of a solution because  $dd^c = -(2\pi i)^{-1}\partial \overline{\partial}.$ If $\mu$ is real-valued, then  the imaginary part of a solution is harmonic. Hence, its real part is a real-valued solution. $\blacksquare$

\bibliographystyle{alpha}

\bibliography{/home/ricardo/Escritorio/Produccion/citas}

\end{document}